\documentclass[11pt]{article} 
\newcommand{\doctype}{TECH}

\usepackage{ifthen}

\newcommand{\Prob}{\ensuremath{\mathbb{P}}}

\long\def\comment#1{}

\usepackage{fullpage}

\newcommand{\citep}[1]{\cite{#1}}

\usepackage{amsfonts,graphicx}
\usepackage{amsfonts,amsmath,amssymb,amsthm}
\usepackage{epic,eepic,eepicemu}
\usepackage{epsf}
\usepackage{graphics}
\usepackage{psfrag}
\usepackage{epsfig}
\usepackage{subfig}

\theoremstyle{plain}
\newtheorem{theorem}{Theorem}

\newtheorem{lemma}{Lemma}

\theoremstyle{definition}
\newtheorem{example}{Example}

\makeatletter
\def\@cite#1#2{[\if@tempswa #2 \fi #1]}
\makeatother

\def\E{{\mathbb E}}
\def\sign{{\rm sign}}

\def\w{{w}}

\def\argmax{{\rm argmax}}
\def\argmin{{\rm argmin}}

\newcommand{\defn}{\ensuremath{:=}}

\newcommand{\real}{\ensuremath{\mathbb{R}}}

\newcommand{\inprod}[2]{\ensuremath{\langle #1 , \, #2 \rangle}}
\newcommand{\convp}[1]{\ensuremath{\stackrel}{#1}{\longrightarrow}}

\newcommand{\convas}{\ensuremath{\stackrel}{a.s.}{\longrightarrow}}

\newcommand{\pnorm}{\ensuremath{r}}

\newcommand{\ftrue}{\ensuremath{f_0}}

\newcommand{\fest}{\ensuremath{\widehat{f}_n}}

\newcommand{\ptrue}{\ensuremath{p_0}}
\newcommand{\qtrue}{\ensuremath{q_0}}
\newcommand{\gtrue}{\ensuremath{g_0}}
\newcommand{\gest}{\ensuremath{\widehat{g}_n}}

\newcommand{\logGenv}{\ensuremath{G_1}}
\newcommand{\Genv}{\ensuremath{G_0}}
\newcommand{\logGnormal}{\ensuremath{G_2}}

\newcommand{\Hspace}{\ensuremath{\mathcal{H}}}
\newcommand{\Xspace}{\ensuremath{\mathcal{X}}}

\newcommand{\ApproxErr}{\ensuremath{\mathcal{E}_0}}
\newcommand{\EstimErr}{\ensuremath{\mathcal{E}_1}}
\newcommand{\Sobspace}[3]{\ensuremath{W_{#1}^{#2}(#3)}}

\newcommand{\dims}{\ensuremath{d}}
\newcommand{\devdeg}{\ensuremath{\alpha}}

\newcommand{\entropy}[3]{\ensuremath{\mathcal{H}_{#1}(#2,#3)}}
\newcommand{\entrobrack}[3]{\ensuremath{\mathcal{H}_{#1}^{B}(#2,#3)}}

\newcommand{\Ker}{\ensuremath{K}}

\newcommand{\vtiny}{\ensuremath{\vspace*{.01in}}}

\newcommand{\myparagraph}[1]{{\noindent \bf{#1}}}

\newcommand{\Pdist}{\ensuremath{\mathbb{P}}}
\newcommand{\Qdist}{\ensuremath{\mathbb{Q}}}

\newcommand{\Gclass}{\ensuremath{\mathcal{G}}}
\newcommand{\Fclass}{\ensuremath{\mathcal{F}}}

\newcommand{\west}{\ensuremath{\widehat{w}}}
\newcommand{\alphaest}{\ensuremath{\widehat{\alpha}}}

\newcommand{\Xsam}[1]{\ensuremath{X^{(#1)}}}
\newcommand{\Ysam}[1]{\ensuremath{Y^{(#1)}}}
\newcommand{\nsam}{\ensuremath{n}}
\newcommand{\msam}{\ensuremath{n}}

\newcommand{\Mfix}{\ensuremath{{M^{*}}}}

\newcommand{\order}{\ensuremath{\mathcal{O}}}

\makeatletter
\long\def\@makecaption#1#2{
        \vskip 0.8ex
        \setbox\@tempboxa\hbox{\small {\bf #1:} #2}
        \parindent 1.5em  
        \dimen0=\hsize
        \advance\dimen0 by -3em
        \ifdim \wd\@tempboxa >\dimen0
                \hbox to \hsize{
                        \parindent 0em
                        \hfil 
                        \parbox{\dimen0}{\def\baselinestretch{0.96}\small
                                {\bf #1.} #2
                                } 
                        \hfil}
        \else \hbox to \hsize{\hfil \box\@tempboxa \hfil}
        \fi
        }
\makeatother

\makeatletter
\long\def\barenote#1{
    \insert\footins{\footnotesize
    \interlinepenalty\interfootnotelinepenalty 
    \splittopskip\footnotesep
    \splitmaxdepth \dp\strutbox \floatingpenalty \@MM
    \hsize\columnwidth \@parboxrestore
    {\rule{\z@}{\footnotesep}\ignorespaces
      #1\strut}}}
\makeatother

\title{{\LARGE{\bf{Estimating divergence functionals and the likelihood
ratio by convex risk minimization}}}}
\author{ XuanLong Nguyen\\ Dept. of Statistical Science \\ 
Duke University\\ \texttt{\small xuanlong.nguyen@stat.duke.edu} \\[.2cm]
\and
Martin J. Wainwright \\ Dept. of Statistics, and Dept. of EECS\\
University of California, Berkeley\\ \texttt{\small
wainwrig@stat.berkeley.edu} \\[.2cm]
\and
Michael I. Jordan \\ Dept. of Statistics, and Dept. of EECS \\
University of California, Berkeley\\ \texttt{\small
jordan@stat.berkeley.edu}\\ }

\vspace{.2cm}

\begin{document}
\date{Revised April 15, 2009}

\maketitle

\ifthenelse{\equal{\doctype}{IEEE}}
{\typeout{}}
{ \begin{center}
Technical Report 764\\
Department of Statistics \\
University of California, Berkeley\\
\bigskip
\end{center}
}

\begin{abstract}

We develop and analyze $M$-estimation methods for divergence
functionals and the likelihood ratios of two probability
distributions. Our method is based on a non-asymptotic variational
characterization of $f$-divergences, which allows the problem of
estimating divergences to be tackled via convex empirical risk 
optimization. The resulting estimators are simple to implement, 
requiring only the solution of standard convex programs.  We present 
an analysis of consistency and convergence for these estimators. 
Given conditions only on the ratios of densities, we show that 
our estimators can achieve optimal minimax rates for the likelihood 
ratio and the divergence functionals in certain regimes.  
We derive an efficient optimization 
algorithm for computing our estimates, and illustrate their convergence 
behavior and practical viability by simulations.\footnote{Preliminary 
versions of this work were presented at the International Symposium on
Information Theory (2007) ~\citep{NguyenEtAl-isit07} and the Neural 
Information Processing Systems Conference (2007)~\citep{NguyenEtAl-nips07}.}

\end{abstract}

\section{Introduction}

Divergences (or pseudo-distances) based on likelihood ratios between
pairs of multivariate probability distribution densities play a
central role in information theory and statistics.  For instance, in
the asymptotic analysis of hypothesis testing, the Kullback-Leibler
and Chernoff divergences control the decay rates of error
probabilities (e.g., see Stein's lemma~\cite{Cover91} and its
variants).  As a particular case of the Kullback-Leibler divergence,
the mutual information specifies capacities in channel coding coding
and data compression~\cite{Cover91}.  In statistical machine learning
and signal processing, divergences between probability distributions
are frequently exploited as metrics to be optimized, such as in
independent component analysis~\cite{Comon94} and decentralized
detection~\cite{Tsitsiklis93}.

In all of these settings, an important problem is that of
\emph{divergence estimation}: how to estimate the divergence between
two multivariate probability distributions, say $\Pdist$ and $\Qdist$,
based on a set of samples from each distribution?  A canonical example
is estimation of the Kullback-Leibler (KL) divergence from samples.
This problem includes as a special case the problem of estimating the
mutual information, corresponding to the KL divergence between a joint
distribution and the product of its marginals, as well as the problem
of estimating the Shannon entropy of a distribution $\Pdist$, which is
related to the KL divergence between $\Pdist$ and the uniform
distribution.  Several researchers have studied the problem of Shannon
entropy estimation~\cite{Gyorfi-Meulen,Joe-1989,Hall-Morton} based on
various types of nonparametric techniques.  Somewhat more generally,
the problem of estimating an integral functional of a single density
has been studied extensively, dating back to early
work~\citep{ibragimov-khasminskii,levit} from 1970s, and continuing on
in later research~\cite{Bickel-Ritov,Birge-Massart,Laurent}.  More
recent recent work by Wang et al.~\cite{Wang05} has developed
algorithms for estimating the KL divergence between a pair of
continuous distributions $\Pdist$ and $\Qdist$, based on building
data-dependent partitions of equivalent (empirical)
$\Qdist$-measure. Wang et al.~\cite{Wang06} also proposed an
interesting nonparametric estimator of the KL divergence using
$1$-nearest neighbor technique. Both estimators were empirically shown
to outperform direct plug-in methods, but no theoretical results on
convergence rates were provided.

In this paper, we propose methods for estimating divergence
functionals as well as likelihood density ratios based on simple
$M$-estimators.  Although our primary interest is the Kullback-Leibler
divergence, our methodology is more broadly applicable to the class of
Ali-Silvey distances, also known as
$f$-divergences~\citep{Ali-Silvey66,Csiszar67}. Any divergence in this
family, to be defined more formally in the sequel, is of the form
$D_\phi(\Pdist,\Qdist) = \int \phi(d \Qdist /d \Pdist) d\Pdist$, where
$\phi$ is a convex function of the likelihood ratio $d\Qdist/
d\Pdist$.

Our estimation method is motivated by a non-asymptotic
characterization of $f$-divergence, due independently to several
authors~\cite{Bron-Keziou,Keziou,NguyenEtAl-AOS}.  Roughly speaking,
the main theorem in~\cite{NguyenEtAl-AOS} states that that
there is a correspondence between the family of $f$-divergences and a
family of losses such that the minimum risk is equal to the negative
of the divergence. In other words, any negative $f$-divergence can
serve as a lower bound for a risk minimization problem.  This
correspondence provides a \emph{variational characterization}, by
which the divergence $D_\phi(\Pdist,\Qdist)$ can be expressed as the
maximum of an Bayes decision problem involving two hypotheses $\Pdist$
and $\Qdist$.  In this way, as we show in this paper, estimating the
divergence $D_\phi(\Pdist,\Qdist)$ has an equivalent reformulation as
solving a certain Bayes decision problem.  This reformulation leads to
an $M$-estimation procedure, in which the divergence is estimated by
solving the convex program defined by the Bayes decision problem.  
This approach not only leads to an $M$-estimation procedure for the
divergence but also for the likelihood ratio $d \Pdist/d\Qdist$.

Our second contribution is to analyze the convergence and consistency
properties of our estimators, under certain assumptions on the
permitted class $\Gclass$ of density ratios, or logarithms of density
ratios.  The analysis makes use of some known results in empirical
process theory for nonparametric density
estimation~\citep{VdG-book,vanderVaart96}.  The key technical
condition is the continuity of the suprema of two empirical processes,
induced by $\Pdist$ and $\Qdist$ distributions respectively, with
respect to a metric defined on the class $\Gclass$ of permitted
functions.  This metric arises as a surrogate lower bound of a Bregman
divergence defined on a pair of density ratios. If $\Gclass$ is 
a smooth function class with smooth parameter $\alpha > \dims/2$,
it can be shown that our estimates of the likelihood ratio and
the KL divergence are both optimal in the minimax sense with 
the rate $n^{-\alpha/(\dims + 2\alpha)}$ and $n^{-1/2}$, respectively.

Our third contribution is to provide an efficient implementation 
of one version of our estimator, in which the function class 
$\Gclass$ is approximated by a reproducing kernel Hilbert 
space (RKHS)~\citep{Saitoh88}.  After computing the convex
dual, the estimator can be implemented by solving a simple convex
program involving only the Gram matrix defined by the kernel 
associated with the RKHS.  Our method thus inherits the simplicity of
other kernel-based methods used in statistical machine
learning~\citep{Scholkopf02,ShaweTaylor-Cristianini}.  We
illustrate the empirical behavior of our estimator on various
instances of KL divergence estimation.

The remainder of this paper is organized as follows. In
Section~\ref{Sec:variational}, we provide the variational
characterization of $f$-divergences in general, and KL divergence in
particular.  We then develop an $M$-estimator for the KL divergence and
the likelihood ratio.  Section~\ref{Sec:consistency} is devoted to the
analysis of consistency and convergence rates of these estimators.  In
Section~\ref{Sec:algorithm}, we develop an efficient kernel-based
method for computing our $M$-estimates, and provide simulation results
demonstrating their performance. In Section~\ref{Sec:general}, we
discuss our estimation method and its analysis in a more general
light, encompassing a broader class of $f$-divergences.  We conclude
in Section~\ref{Sec:conclude}. \\

\vspace*{.04in}

\myparagraph{Notation:} For convenience of the reader, we summarize
some notation to be used throughout the paper.  Given a probability
distribution $\Pdist$ and random variable $f$ measureable with respect
to $\Pdist$, we use $\int f d \Pdist$ to denote the expectation of $f$
under $\Pdist$.  When $\Pdist$ is absolutely continuous with respect
to Lesbesgue measure, say with density $\ptrue$, this integral is
equivalent to the usual Lebesgue integral $\int f \ptrue d \mu = \int
f(x) \ptrue(x) \mu (dx)$.  Given samples $\Xsam{1}, \ldots,
\Xsam{\nsam}$ from $\Pdist$, the empirical distribution $\Pdist_\nsam$
is given by $\Pdist_n = \frac{1}{\nsam} \sum_{i=1}^n
\delta_{\Xsam{i}}$, corresponding to a sum of delta functions centered
at the data points.  We use $\int f d \Pdist_\nsam$ as a convenient
short-hand for the empirical expectation $\frac{1}{\nsam}
\sum_{i=1}^\nsam f(\Xsam{i})$.

\section{$M$-estimators for KL divergence and the density ratio}
\label{Sec:variational}

We begin by defining $f$-divergences, and describing a variational
characterization in terms of a Bayes decision problem.  We then
exploit this variational characterization to develop an $M$-estimator.

\subsection{Variational characterization of $f$-divergence}

Consider two probability distributions $\Pdist$ and $\Qdist$, with
$\Pdist$ absolutely continuous with respect to $\Qdist$.  Assume
moreover that both distributions are absolutely continuous with
respect to Lebesgue measure $\mu$, with densities $\ptrue$ and
$\qtrue$, respectively, on some compact domain $\Xspace \subset
\real^d$.  The Kullback-Leibler (KL) divergence between $\Pdist$ and
$\Qdist$ is defined by the integral
\begin{eqnarray}
\label{EqnDefnKL}
D_{K}(\Pdist,\Qdist) & \defn & \int \ptrue \log \frac{\ptrue}{\qtrue}
\;d\mu.
\end{eqnarray}
This divergence is a special case of a broader class of divergences
known as Ali-Silvey distances or
$f$-divergence~\citep{Csiszar67,Ali-Silvey66}, which take the form
\begin{eqnarray}
\label{EqnDefnFdvi}
D_{\phi}(\Pdist,\Qdist) & \defn & \int \ptrue \phi(\qtrue/\ptrue) \;
d\mu,
\end{eqnarray}
where $\phi: \real \rightarrow \real$ is a convex and lower
semi-continuous (l.s.c.) function.  Different choices of $\phi$ result
in many divergences that play important roles in information theory
and statistics, including not only the KL divergence~\eqref{EqnDefnKL}
but also the total variational distance, the Hellinger distance, and
so on; see~\citep{Topsoe} for further examples.

We begin by stating and proving a variational representation for the
divergence $D_\phi$.  In order to do so, we require some basic
definitions from convex analysis~\cite{Rockafellar,Hiriart1}.  The
\emph{subdifferential} $\partial \phi(t)$ of the convex function $\phi$ at a
point $t \in \real$ is the set
\begin{eqnarray}
\partial \phi(t) & \defn & \left \{ z \in \real \; \mid \; \phi(s) \geq
\phi(t) + z \, (s-t)  \qquad \forall s \in \real \right \}.
\end{eqnarray}
As a special case, if $\phi$ is differentiable at $t$, then $\partial
\phi(t) = \{\phi'(t)\}$.  The function $\phi^*$ is the \emph{conjugate dual
function} associated with $\phi$, defined as
\begin{eqnarray}
\phi^*(v) & \defn & \sup_{u \in \real} \left\{ u \, v - \phi(u) \right \}.
\end{eqnarray}
With these definitions, we have:
\begin{lemma} 
\label{Lem-variational}
For any class of functions $\Fclass$ mapping from $\Xspace$ to
$\real$, we have the lower bound
\begin{equation}
\label{EqnVariational}
D_{\phi}(\Pdist,\Qdist) \geq \sup_{f \in \Fclass} \int \left[f \;
d\Qdist - \phi^{*}(f) \;d\Pdist \right].
\end{equation}
Equality holds if and only if the subdifferential $\partial
\phi(\qtrue/\ptrue)$ contains an element of $\Fclass$.
\end{lemma}
\begin{proof}
Since $\phi$ is convex and lower semi-continuous, Fenchel convex
duality~\citep{Rockafellar} guarantees that we can write $\phi$ in terms
of its conjugate dual as $\phi(u) = \sup_{v\in \real} \big \{uv -
\phi^{*}(v) \big \}$.  Consequently, we have
\begin{eqnarray*}
D_{\phi}(\Pdist,\Qdist) & = & \int \ptrue \sup_f 
(f \qtrue/\ptrue - \phi^{*}(f)) \; d\mu \\ 
& = & \sup_{f} \int \left [f\qtrue - \phi^{*}(f)\ptrue \right] \;d\mu
\\ 
& = & \sup_{f} \int \left [f \; d\Qdist - \phi^{*}(f) \;d\Pdist \right],
\end{eqnarray*}
where the supremum is taken over all measurable functions $f:\Xspace
\rightarrow \real$.  It is simple to see that equality in the 
supremum is attained at a function $f$ such that $\qtrue/\ptrue \in
\partial \phi^{*}(f)$ where $\qtrue,\ptrue$ and $f$ are evaluated at
any $x \in \Xspace$. By convex duality, this is true if $f \in
\partial \phi(\qtrue/\ptrue)$ for any $x\in \Xspace$. 
\end{proof}

\vspace*{.1in}

\subsection{$M$-estimators for the KL divergence and likelihood ratio}

We now describe how the variational
representation~\eqref{EqnVariational} specializes to an $M$-estimator
for the Kullback-Leibler (KL) divergence.  As a particular
$f$-divergence, the KL divergence is induced by the convex function
\begin{eqnarray}
\phi(t) & = & \begin{cases} - \log(t) & \mbox{for $t > 0$} \\
                            +\infty & \mbox{for $t \leq 0$}.
	      \end{cases}
\end{eqnarray}
A short calculation shows that the conjugate dual takes the form
\begin{eqnarray}
\phi^*(v) & = & \begin{cases}  -1 - \log(-v) & \mbox{if $u< 0$, and} \\
+\infty & \mbox{otherwise.}
		\end{cases}
\end{eqnarray}
As a consequence of Lemma~\ref{Lem-variational}, we obtain the
following representation of the Kullback-Leibler divergence:
$D_K(\Pdist,\Qdist) = \sup_{f<0} \int f\; d\Qdist - \int [-1
-\log(-f)] d\Pdist$.  
After the substitution $g = -f$, this can be written as
\begin{equation}
\label{Eqn-KL}
D_K(\Pdist, \Qdist) = \sup_{g>0} \int \log g \; d\Pdist - \int g
d\Qdist + 1,
\end{equation}
for which the supremum is attained at $\gtrue = \ptrue/\qtrue$.  

We now take a statistical perspective on the variational 
problem~\eqref{Eqn-KL}, where we assume that the
distributions $\Pdist$ and $\Qdist$ are unknown.  
We suppose that we are given independent and identically distributed
(i.i.d.)  samples, say $\Xsam{1}, \Xsam{2}, \ldots, \Xsam{\nsam}$
drawn i.i.d.\ from $\Pdist$, and $\Ysam{1}, \Ysam{2}, \ldots,
\Ysam{\msam}$ drawn i.i.d.\ from $\Qdist$.  Denote by $\Pdist_\nsam$
the empirical distribution defined by the samples $\{\Xsam{1}, \ldots,
\Xsam{\nsam} \}$, given explicitly by $\Pdist_\nsam \defn
\frac{1}{\nsam} \sum_{i=1}^\nsam \delta_{\Xsam{i}}$, with the
empirical distribution $\Qdist_\msam$ associated with $\{\Ysam{1},
\ldots, \Ysam{\msam} \}$ defined analogously.  We consider
two classes of estimators:

\paragraph{Estimator E1:}
Given the empirical distributions, we consider the estimator 
obtained by replacing the true distributions $\Pdist$ and $\Qdist$ 
with their empirical versions, and maximizing over some pre-specified
class $\Gclass$ of functions $g:\Xspace \rightarrow \real_+$, as follows:
\begin{subequations}
\begin{eqnarray}
\label{Eqn-g}
\gest = \argmax_{g\in \Gclass} \int \log g \; d\Pdist_n - \int g \; 
d\Qdist_n,\\
\label{Eqn-KLest}
\widehat{D}_K = \int \log \gest \; d\Pdist_n - \int \gest \;
d\Qdist_n + 1.
\end{eqnarray}
\end{subequations}
Assuming that $\Gclass$ is a convex set of functions, the 
implementation of the estimator requires solving 
a convex optimization problem, albeit over an
(infinite-dimensional) function space. For this estimation method to
work, several conditions on $\Gclass$ are required. First, so as to
control approximation error, it is natural to require that $\Gclass$
is sufficiently rich so as to contain the true likelihood ratio
$\gtrue$ in the sense of KL divergence, i.e., there is some $g \in
\Gclass$ such that $g = \gtrue$ a.e..  On the other hand, $\Gclass$
should not be too large, so that estimation is possible.  To formalize
this condition, let $I(g)$ be a measure of complexity for $g$, where
$I$ is a non-negative functional and $I(\gtrue) < \infty$.  Given some
fixed finite constant $\Mfix \geq I(\gtrue)$, we then define
\begin{equation}
\label{Eqn-Gclass-E1}
\Gclass \defn \Gclass_\Mfix \defn \{g: I(g) \leq \Mfix\}.
\end{equation}

\vspace*{.1in}

\paragraph{Estimator E2:}
In practice, the ``true'' $I(\gtrue)$ is not known, and hence it is
not clear as a practical matter how to choose the fixed $\Mfix$ 
defining estimator E1.  Thus we also consider an approach that involves 
an explicit penalty $I(g)$. In this approach, let
\begin{equation}
\label{EqnGclassDecomp}
\Gclass = \cup_{1\leq M \leq \infty} \Gclass_M.
\end{equation}
The estimation procedure involves solving the following program:
\begin{subequations}
\begin{eqnarray}
\label{Eqn-penalize}
\gest & = & \argmin_{g\in \Gclass} \int g d\Qdist_n - \int \log g \;
d\Pdist_n + \frac{\lambda_n}{2}I^2(g),\\
\label{Eqn-penalize-KLest}
\widehat{D}_K & = & \int \log \gest \; d\Pdist_n - \int \gest \;
d\Qdist_n + 1,
\end{eqnarray}
\end{subequations}
where $\lambda_n > 0$ is a regularization parameter.   \\

\vspace*{.1in}

As we discuss in Section~\ref{Sec:algorithm}, for function classes
defined by reproducing kernel Hilbert spaces,
problems~\eqref{Eqn-KLest} and~\eqref{Eqn-penalize} can actually be
posed as a finite-dimensional convex programs (in $\nsam$ dimensions),
and solved efficiently by standard methods.  In addition to the
estimate $\widehat{D}_K$ of the KL divergence, if the
supremum is attained at $\gest$, then $\gest$ is an $M$-estimator of the
density ratio $\gtrue \defn \ptrue/\qtrue$.

In the next section, we present results regarding the consistency and
convergence rates of both estimators. While these methods have similar
convergence behavior, estimator E1 is somewhat simpler to analyze and admits
weaker conditions for consistency. On the other hands, estimator E2
seems more practical. Details of algorithmic derivations for estimator
E2 are described in Section~\ref{Sec:algorithm}.

\section{Consistency and convergence rate analysis}
\label{Sec:consistency}

For the KL divergence functional, the difference $|\widehat{D}_K -
D_K(\Pdist,\Qdist)|$ is a natural performance measure. For estimating
the density ratio function, this difference can also be used,
although more direct metrics are customarily preferred.
In our analysis, we view $\gtrue = \ptrue/\qtrue$ as a density function with
respect to $\Qdist$ measure, and adopt the (generalized) Hellinger
distance as a performance measure for estimating the likelihood ratio
function:
\begin{equation}
h_\Qdist^2(\gtrue, g) \defn \frac{1}{2} \int (\sqrt{\gtrue} -
\sqrt{g})^2 \; d\Qdist.
\end{equation}

\subsection{Consistency of estimator E1}

Our analysis of consistency relies on tools from empirical process 
theory.  Let us briefly review the notion of the \emph{metric entropy} of 
function classes (see, e.g.,~\citep{vanderVaart96} for further background).
For any $r \geq 1$ and distribution function $\Qdist$, define the
empirical $L_r$ metric
\begin{eqnarray*}
\| g \|^r_{L_r(\Qdist)} & \defn & \int |g|^r d\Qdist,
\end{eqnarray*}
and let $L_r(\Qdist)$ denote the metric space defined by this
distance.  For any fixed $\delta > 0$, a covering for function class
$\Gclass$ using the metric $L_r(\Qdist)$ is a collection of functions
which allow $\Gclass$ to be covered using $L_r(\Qdist)$ balls of radius
$\delta$ centered at these functions.  Letting $N_\delta(\Gclass,
L_r(\Qdist))$ be the smallest cardinality of such a covering, then
$\entropy{\delta}{\Gclass}{L_r(\Qdist)} \defn \log
N_\delta(\Gclass,L_r(\Qdist))$ is called the entropy for $\Gclass$
using the $L_r(\Qdist)$ metric. A related notion is \emph{entropy with
bracketing}.  Let $N_\delta^B(\Gclass,L_r(\Qdist))$ be the smallest
value of $N$ for which there exist pairs of functions
$\{g_j^L,g_j^U\}$ such that $\|g_j^U - g_j^L\|_{L_r(\Qdist)} \leq
\delta$, and such that for each $g\in \Gclass$ there is a $j$ such that
$g_j^L \leq g \leq g_j^L$. Then
$\entrobrack{\delta}{\Gclass}{L_r(\Qdist)} \defn \log
N_\delta^B(\Gclass,L_r(\Qdist))$ is called the entropy with bracketing
of $\Gclass$.
Define the envelope functions:
\begin{equation}
 \Genv(x) = \sup_{g \in \Gclass} |g(x)|, \qquad \mbox{and} \qquad
 \logGenv(x) = \sup_{g \in \Gclass} |\log\frac{g(x)}{\gtrue(x)}|.
\end{equation}

For the estimator E1, we impose the following assumptions on the
distributions $\Pdist,\Qdist$ and the function class $\Gclass$.
\paragraph{Assumptions.} 
\begin{enumerate}
\item[(A).]
The KL divergence is bounded: $D_K(\Pdist,\Qdist) < \infty$.
\item[(B).] There is some $g \in \Gclass$
such that $g = \gtrue$ almost surely (a.s.).
\end{enumerate}

In the following theorem, the almost sure statement can be taken with
respect to either $\Pdist$ or $\Qdist$ since they share the same
support.

\begin{theorem} 
\label{ThmEstimE1}
Suppose that assumptions (A) and (B) hold. \\
\noindent (a) Assume the envelope conditions 
\begin{subequations}
\label{Eqn-envelope}
\begin{eqnarray}
\label{Eqn-envelope-Q}
\int \Genv d\Qdist & < & \infty \\
\label{Eqn-envelope-P}
\int \logGenv d\Pdist & < & \infty 
\end{eqnarray}
\end{subequations} 
and suppose that for all $\delta > 0$ there holds:
\begin{subequations}
\label{Eqn-entropy}
\begin{eqnarray}
\label{Eqn-entropy-Q}
\frac{1}{n} \entropy{\delta}{\Gclass - \gtrue}{L_1(\Qdist_n)} 
\convp{\Qdist} 0, && \\
\label{Eqn-entropy-P}
\frac{1}{n} \entropy{\delta}{\log \Gclass/\gtrue}{L_1(\Pdist_n)} 
\convp{\Pdist} 0.&& 
\end{eqnarray}
\end{subequations}
Then, $|\widehat{D}_K - D_K(\Pdist,\Qdist)| \convas 0$,
and $h_\Qdist(\gtrue,\gest) \convas 0$. \\

\noindent (b)  Suppose only that~\eqref{Eqn-envelope-Q}
and~\eqref{Eqn-entropy-Q} hold, and
\begin{equation}
\label{Eqn-entropy-P-2}
\frac{1}{n} \entropy{\delta}{\log \frac{\Gclass + \gtrue}{2\gtrue}}{L_1(\Pdist_n)} 
\convp{\Pdist} 0.
\end{equation}
Then $h_\Qdist(\gtrue,\gest) \convas 0$.
\end{theorem}

To provide intuition for the conditions of Theorem~\ref{ThmEstimE1},
note that the envelope condition~\eqref{Eqn-envelope-Q} is relatively
mild, satisfied (for instance) if $\Gclass$ is uniformly bounded from
above. On the other hand, the envelope
condition~\eqref{Eqn-envelope-P} is much more stringent.  Due to
the logarithm, this can be satisfied if all functions in $\Gclass$ are
bounded from \emph{both} above and below.  However, as we see
in part (b), we do not require boundedness from below; to ensure
Hellinger consistency we can drop both the envelope 
condition~\eqref{Eqn-envelope-P} and the entropy
condition~\eqref{Eqn-entropy-P}, replacing them with the milder
entropy condition~\eqref{Eqn-entropy-P-2}.

It is worth noting that both~\eqref{Eqn-entropy-Q}
and~\eqref{Eqn-entropy-P-2} can be deduced from the following single
condition: for all $\delta > 0$, the bracketing entropy satisfies
\begin{equation}
\label{Eqn-entropy-bracket}
\entrobrack{\delta}{\Gclass}{L_1(\Qdist)} < \infty.
\end{equation}
Indeed, given equation~\eqref{Eqn-envelope-Q} and by the law of large
numbers, condition~\eqref{Eqn-entropy-bracket} directly
implies~\eqref{Eqn-entropy-Q}.  To establish
condition~\eqref{Eqn-entropy-P-2}, note that by a Taylor expansion,
we have
\begin{equation*}\biggr |\log \frac{g_1 + \gtrue}{2\gtrue} - \log \frac{g_2 + \gtrue}{2\gtrue}\biggr |
\; \leq \; \frac{|g_1 - g_2|}{\gtrue},
\end{equation*}
so that $\frac{1}{n}\entropy{\delta}{\log \frac{\Gclass +
\gtrue}{2\gtrue}}{L_1(\Pdist_n)} \leq
\frac{1}{n}\entropy{\delta}{\Gclass/\gtrue}{L_1(\Pdist_n)}$.  Since
$\Genv \in L_1(\Qdist)$, we have $\Genv/\gtrue \in L_1(\Pdist)$. In
addition, $\entrobrack{\delta}{\Gclass/\gtrue}{L_1(\Pdist)} \leq
\entrobrack{\delta}{\Gclass}{L_1(\Qdist)} < \infty$. By the law of
large numbers, the metric entropy
$\entropy{\delta}{\Gclass/\gtrue}{L_1(\Pdist_n)}$ is bounded in
probability, so that condition~\eqref{Eqn-entropy-P-2} holds.

In practice, the entropy conditions are satisfied by a variety of
function classes. Examples include various types of reproducing
kernel Hilbert spaces~\cite{Saitoh88}, as described in more detail
in Section~\ref{Sec:algorithm}, as well as the Sobolev classes,
which we describe in the following example.

\begin{example}[Sobolev classes]
\label{ExaSobolev}
Let $\kappa = (\kappa_1,\ldots,\kappa_d)$ be a $d$-dimensional
multi-index, where all $\kappa_i$ are natural numbers.  Given a vector
$x \in \real^d$, define $x^\kappa \; \defn \;
\prod_{i=1}^{d}x_i^{\kappa_i}$ and $|\kappa| \defn
\sum_{i=1}^{d}\kappa_i$. For a suitably differentiable function $f$,
let $D^\kappa$ denote the multivariate differential operator
\begin{equation}
D^\kappa f(x) \; \defn \; \frac{\partial^{|\kappa|}}{\partial
x_1^{\kappa_1} \ldots \partial x_d^{\kappa_d}} f(x_1,\ldots,x_d),
\end{equation}
and define the norm $\|f\|_{L_2^{\devdeg}(\Xspace)}^{2} \; \defn \;
\sum_{|\kappa| = \devdeg} \int |D^{\kappa} f(x)|^2 \; dx$.

With this notation, we define the norm
\begin{equation}
\|f\|_{\Sobspace{2}{\devdeg}{\Xspace}} \; \defn \; \|f
\|_{L_2(\Xspace)} + \|f \|_{L_2^{\devdeg}(\Xspace)},
\end{equation}
and the \emph{Sobolev space} $\Sobspace{2}{\devdeg}{\Xspace}$ of
functions with finite $\|f\|_{\Sobspace{2}{\devdeg}{\Xspace}}$-norm.
Suppose that the domain $\Xspace$ is a compact subset of
$\real^\dims$.  Let the complexity measure $I$ be the Sobolev
norm---that is, \mbox{$I(g) \; \defn \;
\|g\|_{\Sobspace{\pnorm}{\devdeg}{\Xspace}}$.}  With this choice of
complexity measure, it is known~\cite{Birman-Solomjak} that the
function class $\Gclass$ defined in equation~\eqref{Eqn-Gclass-E1}
satisfies, for any $\delta > 0$, the metric entropy bound
\begin{eqnarray}
\entropy{\delta}{\Gclass}{L_\infty} & = & \order
\left(\delta^{-\dims/\alpha} \right) < \infty,
\end{eqnarray}
for all smoothness indices $\alpha > \dims/2$.  As a result, both
conditions~\eqref{Eqn-entropy-bracket} and~\eqref{Eqn-entropy-P} hold
if, for instance, $\Gclass$ is restricted to a subset of smooth
functions that are bounded from above, and $\gtrue$ is bounded from
below. %
\hfill $\diamondsuit$
\end{example}

\subsection{Proof of Theorem~\ref{ThmEstimE1}}

We now turn to the proof of Theorem~\ref{ThmEstimE1}, beginning with
part (a).  Define the following quantities associated with the
function class $\Gclass$:
\begin{eqnarray}
\label{Eqn-ApproxErr} 
\ApproxErr(\Gclass) & = & D_K(\Pdist,\Qdist) - \sup_{g\in \Gclass}
\int (\log g \; d\Pdist - g \;d\Qdist + 1) \geq 0 \\
\label{Eqn-EstimErr}
\EstimErr(\Gclass) & = & \sup_{g \in \Gclass} \biggr |\int 
\log g \; d(\Pdist_n - \Pdist)\;  - g d(\Qdist_n - \Qdist) \biggr |.
\end{eqnarray}
The quantity $\ApproxErr$ is the \emph{approximation error}, which
measures the bias incurred by limiting the optimization to the class
$\Gclass$.  The term $\EstimErr$ is the \emph{estimation error}
associated with the class $\Gclass$.  Our focus in this paper
is the statistical problem associated with the estimation error
$\EstimErr$, and thus we have imposed assumption (B), which implies 
that the approximation error $\ApproxErr(\Gclass) = 0$.  
Moreover, from equations~\eqref{Eqn-KL} and~\eqref{Eqn-KLest}, 
straightforward algebra yields that
\begin{eqnarray}
\label{EqnBound}
|\widehat{D}_K - D_K(\Pdist, \Qdist)| & \leq & \ApproxErr(\Gclass) +
 \EstimErr(\Gclass) =  \EstimErr(\Gclass).
\end{eqnarray}
Accordingly, the remainder of the proof is devoted to analysis of the
estimation error $\EstimErr(\Gclass)$.

In order to analyze $\EstimErr(\Gclass)$, define the following
processes:
\begin{subequations}
\begin{eqnarray}
\label{Eqn-supreme}
v_n(\Gclass) & \defn & \sup_{g\in \Gclass} \biggr |\int \log
\frac{g}{\gtrue} d(\Pdist_n - \Pdist) - \int (g - \gtrue) d(\Qdist_n -
\Qdist) \biggr |, \qquad \mbox{and} \\
w_n(\gtrue) & \defn & \biggr |\int \log \gtrue \; d(\Pdist_n -
\Pdist)\; - \gtrue d(\Qdist_n - \Qdist) \biggr |.
\end{eqnarray}
\end{subequations}
Note that we have
\begin{equation}
\label{Eqn-EstBound}
\EstimErr(\Gclass) \leq v_n(\Gclass) + w_n(\gtrue).
\end{equation}

Note that the quantity $w_n(\gtrue)$ is the difference
between an empirical and a population expectation.  Let us verify that
the conditions for the strong law of large numbers (SLN) apply.  Using
the inequality
\begin{equation*}
\int \ptrue|\log(\ptrue/\qtrue)| \leq D_K(\Pdist,\Qdist) + 4
\sqrt{D_K(\Pdist,\Qdist)},
\end{equation*}
due to Csisz\'ar (cf.~\cite{Gyorfi-Meulen}), it follows that $\log
\gtrue$ is $\Pdist$ integrable. In addition, the function $\gtrue$ is
$\Qdist$ integrable, since $\int \gtrue d \Qdist = \int
(\ptrue/\qtrue) d\Qdist = 1$.  Thus, the SLN applies, and we conclude
that $w_n(\gtrue) \convas 0$.  By applying
Theorem~\ref{Thm-VdG-consistency} from Appendix~\ref{AppEmpProcess},
we conclude that $v_n(\Gclass) \convas 0$.  From the decomposition in
equation~\eqref{Eqn-EstBound}, we conclude that $\EstimErr(\Gclass)
\convas 0$, so that $|\widehat{D}_K - D_K(\Pdist, \Qdist)| \convas 0$.

To establish Hellinger consistency of the likelihood ratio, we require
the following lemma, whose proof is in Appendix~\ref{AppLemDistKL}:
\begin{lemma}
\label{Prop-dist-KL}
Defining the ``distance'' $d(\gtrue,g) \, \defn \, \int (g-\gtrue)
d\Qdist - \log\frac{g}{\gtrue} d\Pdist$, the following statements
hold:
\begin{enumerate}
\item[(i)] For any $g\in \Gclass$, we have $d(\gtrue,g) \geq 2
h_\Qdist^2(g,\gtrue).$
\item[(ii)] For the estimate $\gest$ defined in
equation~\eqref{Eqn-g}, we have $d(\gtrue,\gest) \leq v_n(\Gclass)$.
\end{enumerate}
\end{lemma}
\noindent The Hellinger consistency $h_\Qdist(\gtrue,\gest)
\convas 0$ of Theorem~\ref{ThmEstimE1}(a) is an immediate consequence
of this lemma. \\

Turning now to Theorem~\ref{ThmEstimE1} (b), we require a more refined
lemma relating the Hellinger distance to suprema of empirical
processes.
\begin{lemma} 
\label{Lem-hellinger}
If $\gest$ is an estimate of $g$, then:
\[\frac{1}{8} h_\Qdist^2(\gtrue,\gest) \leq 2h_\Qdist^2(\gtrue, 
\frac{\gtrue + \gest}{2}) \leq - \int \frac{\gest - \gtrue}{2}
d(\Qdist_n - \Qdist) + \int \log \frac{\gest +\gtrue}{2\gtrue}
d(\Pdist_n - \Pdist).\]
\comment{ (ii) If $\Gclass$ is a convex class, and $\gest$ is an
estimate of $g$, then:
\[h_\Qdist^2(\gest, \gtrue) \leq 
-\int \frac{\gest - \gtrue}{2} d(\Qdist_n - \Qdist) + \int
\frac{2\gest}{\gest + \gtrue} d(\Pdist_n - \Pdist).\] }
\end{lemma}
See Appendix~\ref{AppLemHellinger} for the proof of this claim.  To
complete the proof of Theorem~\ref{ThmEstimE1}, define $\logGnormal(x)
= \sup_{g \in \Gclass} |\log\frac{g(x)+\gtrue(x)}{2\gtrue(x)}|$.
Due to Lemma~\ref{Lem-hellinger} and standard results from empirical 
process theory (see Appendix~\ref{AppEmpProcess}, 
Theorem~\ref{Thm-VdG-consistency}) it is sufficient to 
prove that $\int \logGnormal d\Pdist < \infty$.
To establish this claim, note that
\begin{eqnarray*}
\int \logGnormal d\Pdist & \leq & \int \sup_{g\in \Gclass} \max \left
\{ \frac{g(x) + \gtrue(x)}{2\gtrue(x)} -1, \; \log 2 \right\} d\Pdist
\\
& \leq & \log 2 + \int \sup_{g\in \Gclass} |g(x) - \gtrue(x)| d\Qdist
\; \stackrel{(a)}{<} \; \infty,
\end{eqnarray*}
where the last inequality (a) is due to envelope
condition~\eqref{Eqn-envelope-Q}.

\subsection{Convergence rates}
\label{Sec:rates}

In this section, we describe convergence rates for both estimator E1
and estimator E2.  The key technical tool that we use to analyze 
the convergence rate for the likelihood ratio estimate is
Lemma~\ref{Lem-hellinger}, used previously in the proof of
Theorem~\ref{ThmEstimE1}.  This lemma bounds the Hellinger distance
$h_\Qdist(\gtrue,\gest)$ in terms of the suprema of two empirical
processes writh respect to $\Pdist$ and $\Qdist$. In a nutshell, the
suprema of these two empirical processes can be bounded from above in
terms of the Hellinger distance, which allows us to obtain the rates at
which the Hellinger distance goes to zero.

\subsubsection{Convergence rates for estimator E1}
In order to characterize convergence rates for the estimator E1,
we require one of the following two conditions:
\begin{subequations}
\begin{eqnarray}
\label{Eqn-Upbound}
\sup_{g \in \Gclass} \|g\|_\infty \; < \; K_2 && \\
\label{Eqn-UpBelBound}
0 \, < \, K_1 \, \leq \, \inf_x g(x), \quad \mbox{and} \quad \sup_{x}
g(x) \leq K_2 \quad \mbox{for all} \quad g\in \Gclass. &&
\end{eqnarray}
\end{subequations}
We also require the following assumption on function class
$\bar{\Gclass} \defn \{((g+\gtrue)/2)^{1/2}, g\in \Gclass\}$: for some
constant $0 < \gamma_{\bar{\Gclass}} < 2$, there holds for any $\delta
> 0$,
\begin{equation}
\label{Eqn-bounded-entropy-3}
\entrobrack{\delta}{\bar{\Gclass}}{L_2(\Qdist)}
= O(\delta^{-\gamma_{\bar{\Gclass}}}).
\end{equation}
Combining this metric entropy decay rate with
condition~\eqref{Eqn-Upbound}, we deduce that for $\Gclass$, the
bracketing entropy satisfies
\begin{equation}
\entrobrack{\delta}{\Gclass}{L_2(\Qdist)} \; = \;
\order(\delta^{-\gamma_{\bar{\Gclass}}}).
\end{equation}
With these definitions, we can now state a result characterizing the
convergence rate of estimator E1, where the notation $O_\Pdist$ means
``bounded in probability'' with respect to $\Pdist$ measure.
\begin{theorem}[Convergence rates for estimator E1] 
\label{ThmConvE1}
\begin{enumerate}
\item[(a)] If conditions~\eqref{Eqn-Upbound}
and~\eqref{Eqn-bounded-entropy-3} hold, then $h_\Qdist(\gtrue,\gest) =
O_\Pdist(n^{-1/(\gamma_{\bar{\Gclass}} +2)})$.
\item[(b)] If conditions~\eqref{Eqn-UpBelBound}
and~\eqref{Eqn-bounded-entropy-3} hold, then $|\widehat{D}_K -
D_K(\Pdist,\Qdist)| = O_\Pdist(n^{-1/2})$.
%
%
\end{enumerate}
\end{theorem}

\vtiny

\myparagraph{Remarks:} In order to gain intuition for the convergence
rate in part (a), it can be instructive to compare to the
\emph{minimax rate}
\begin{equation*}
r_n \; \defn \; \inf_{\gest \in \Gclass} \sup_{\Pdist,\Qdist}
\E_\Pdist\, [h_\Qdist(\gtrue, \gest)],
\end{equation*}
where the supremum is taken over all pairs $(\Pdist, \Qdist)$ such
that $\gtrue \in \Gclass$.
As a concrete example, if we take $\Gclass$ as the 
Sobolev family from Example~\ref{ExaSobolev}, and if
condition~\eqref{Eqn-UpBelBound} holds, then the minimax rate is $r_n
= \Omega(n^{-1/(\gamma+2)})$, where $\gamma = \gamma_{\bar{\Gclass}} =
\dims/\alpha$ (see  Appendix~\ref{AppMiniMax}).  
Thus, we see that for the Sobolev classes, the
estimator E1 achieves the minimax rate for estimating the likelihood
ratio in Hellinger distance. In addition, the rate $n^{-1/2}$
for the divergence estimate is also optimal.


\subsubsection{Convergence rates for estimator E2}

We now turn to a discussion of the convergence rates of estimator E2.
To analyze this estimator, we assume that
\begin{equation}
\label{Ass-finite-I}
I(\gtrue) < \infty,
\end{equation}
and moreover we assume that the true likelihood ratio $\gtrue$---but \emph{not}
necessarily all of $\Gclass$---is bounded from above and below:
\begin{equation}
\label{Ass-lower-bound}
0 < \eta_0 \leq \gtrue \leq \eta_1 \; \mbox{for some constants}\;
\eta_0,\eta_1.
\end{equation}
We also assume that the sup-norm over $\Gclass_M$ is Lipschitz with
respect to the penalty measure $I(g)$, meaning that there is a constant $c < \infty$ such that for each $M \geq 1$, we have
\begin{equation}
\label{Ass-supnorm-penalize}
\sup_{g\in \Gclass_M} \|g\|_{\infty} \leq c M.
\end{equation}
Finally, we assume that the bracketing entropy of $\Gclass$ satisfies,
for some $0<\gamma < 2$, the condition
\begin{equation}
\label{Ass-entropy}
\entrobrack{\delta}{\Gclass_M}{L_2(\Qdist)} =
\order\left[(M/\delta)^{\gamma} \right] \quad \mbox{for any} \quad
\delta > 0.
\end{equation}
Given these assumptions, we then state the following convergence rate
result for the estimator E2:
\begin{theorem} 
\label{ThmConvE2}
(a) Suppose that assumptions~\eqref{Ass-finite-I}
through~\eqref{Ass-entropy} hold, and that the regularization
parameter $\lambda_n \rightarrow 0$ is chosen such that
\begin{equation*}
\lambda_n^{-1} = \order_\Pdist \left(n^{2/(2+\gamma)} \right).
\end{equation*}
Then under $\Pdist$, we have
\begin{equation}
h_\Qdist(\gtrue,\gest) = \order_\Pdist(\lambda_n^{1/2}).
\end{equation}
(b) Suppose that in addition to assumptions~\eqref{Ass-finite-I}
through~\eqref{Ass-entropy}, there holds \mbox{$\inf \limits_{x \in
\Xspace} \inf \limits_{g\in \Gclass} g(x) \geq \eta_0$.}  Then we have
\begin{equation}
|\widehat{D}_K - D_K(\Pdist,\Qdist)| = \order_\Pdist(n^{-1/2}).
\end{equation}
\end{theorem}

\vspace*{.1in}

\myparagraph{Remarks:} Note that with the choice $\lambda_n^{-1} =
\order_\Pdist \left(n^{2/(2+\gamma)} \right)$ and the special case of
$\Gclass$ as the Sobolev space $\Sobspace{2}{\devdeg}{\Xspace}$ (see
Example~\ref{ExaSobolev}), estimator E2 again achieves the minimax
rate for estimating the density ratio in Hellinger distance.

\subsection{Proof of convergence theorems}

In this section we present a proof of Theorem~\ref{ThmConvE2}. The
proof of Theorem~\ref{ThmConvE1} is similar in spirit, and is provided
in Appendix~\ref{AppThmConvE1}.  The key to our analysis of the
convergence rate of estimator E2 is the following lemma, which can be
viewed as the counterpart of Lemma~\ref{Lem-hellinger}:
\begin{lemma} 
\label{Lem-basic-penalize}
If $\gest$ is an estimate of $g$ using~\eqref{Eqn-penalize}, then:
\begin{equation*}
\frac{1}{4}h_\Qdist^2(\gtrue, \gest) + \frac{\lambda_n}{2}I^2(\gest)
\leq -\int (\gest - \gtrue) d(\Qdist_n - \Qdist) + \int 2\log
\frac{\gest +\gtrue}{2\gtrue} d(\Pdist_n - \Pdist) +
\frac{\lambda_n}{2}I^2(\gtrue) .
\end{equation*}
\end{lemma}
\noindent See Appendix~\ref{AppLemBasicPenalize} for the proof of this
lemma.  Equipped with this auxiliary result, we can now prove
Theorem~\ref{ThmConvE2}(a).  Define $f_g \defn \log
\frac{g+\gtrue}{2\gtrue}$, and let $\Fclass_M \defn \{f_g \; \mid \;
g\in \Gclass_M\}$.  Since $f_g$ is a Lipschitz function of $g$,
conditions~\eqref{Ass-lower-bound} and~\eqref{Ass-entropy} imply that
\begin{equation}
\label{Ass-entropy-penalize}
\entrobrack{\delta}{\Fclass_M}{L_2(\Pdist)} = \order \left \{
(M/\delta)^{\gamma} \right \}.
\end{equation}
Applying Lemma 5.14 from van de Geer~\cite{VdG-book} using the
distance $d_2(\gtrue,g) = \|g-\gtrue\|_{L_2(\Qdist)}$, we have that
the following statement holds under $\Qdist$, and hence holds under
$\Pdist$ as well, since $d \Pdist/d\Qdist$ is bounded from above:
\begin{eqnarray}
\label{Eqn-empiric-3}
\sup_{g \in \Gclass} \frac{ |\int (g-\gtrue) d(\Qdist_n - \Qdist)|}
{n^{-1/2}{d_2(\gtrue,g)}^{1-\gamma/2}(1+I(g) + I(\gtrue))^{\gamma/2}
\vee n^{-\frac{2}{2+\gamma}}(1+I(g) + I(\gtrue)) } =
\order_\Pdist(1).&&
\end{eqnarray} 
In the same vein, we obtain that under $\Pdist$,
\begin{eqnarray}
\label{Eqn-empiric-1}
\sup_{g \in \Gclass} \frac{| \int f_g d(\Pdist_n - \Pdist)|}
{n^{-1/2}{d_2(\gtrue,g)}^{1-\gamma/2}(1+I(g) + I(\gtrue))^{\gamma/2}
\vee n^{-\frac{2}{2+\gamma}}(1+I(g) + I(\gtrue)) } =
\order_\Pdist(1). &&
\end{eqnarray}

Now using condition~\eqref{Ass-supnorm-penalize},  it can be verified
that
\begin{equation*}
d_2(\gtrue,g) = \|g-\gtrue\|_{L_2(\Qdist)} \leq 2c^{1/2}(1 + I(g) +
I(\gtrue))^{1/2} h_\Qdist(\gtrue,g).
\end{equation*}
Combining Lemma~\ref{Lem-basic-penalize} and
equations~\eqref{Eqn-empiric-1} and~\eqref{Eqn-empiric-3}, we conclude
that
\begin{multline}
\label{Eqn-nicebound}
\frac{1}{4}h_\Qdist^2(\gtrue, \gest) + \frac{\lambda_n}{2}I^2(\gest)
\leq \lambda_n I(\gtrue)^2/2 \\ + O_\Pdist \biggr(
n^{-1/2}h_\Qdist(\gtrue,g)^{1-\gamma/2} (1 + I(g) + I(\gtrue))^{1/2 +
\gamma/4} \; \vee \; n^{-\frac{2}{2+\gamma}}(1+I(g) + I(\gtrue))
\biggr ).
\end{multline}
From this point, the proof involves simple algebraic manipulation
of~\eqref{Eqn-nicebound}.
To simplify notation, let $\widehat{h} = h_\Qdist(\gtrue,\gest)$,
$\widehat{I} = I(\gest)$, and $I_0 = I(\gtrue)$.  We break the
analysis into four cases, depending on the behavior of $\widehat{h}$
and $\widehat{I}$: \\

\myparagraph{Case A:}.  In this case, we assume $\widehat{h} \geq
n^{-1/(2+\gamma)} (1 + \widehat{I} + I_0)^{1/2}$ and $\widehat{I} \geq
1 + I_0$.  From~\eqref{Eqn-nicebound}, we have either
\begin{equation*}
\widehat{h}^2/4 + \lambda_n \widehat{I}^2/2 \leq O_\Pdist(n^{-1/2})
\widehat{h}^{1-\gamma/2}\widehat{I}^{1/2 + \gamma/4} \;\; \mbox{or}
\;\; \widehat{h}^2/4 + \lambda_n \widehat{I}^2/2 \leq \lambda_n
I_0^2/2.
\end{equation*}
These conditions imply, respectively, either
\begin{equation*}
\widehat{h} \leq \lambda_n^{-1/2}O_\Pdist(n^{-2/(2+\gamma)}), \;\;\;
\widehat{I} \leq \lambda_n^{-1} O_\Pdist(n^{-2/(2+\gamma)}),
\end{equation*}
or
\begin{equation*}
\widehat{h}\leq O_\Pdist(\lambda_n^{1/2}I_0), \;\;\; \widehat{I} 
\leq O_\Pdist(I_0).
\end{equation*}
In either case, we conclude the proof by setting $\lambda_n^{-1} =
\order_\Pdist(n^{2/(\gamma+2)}(1+I_0))$. \\

\myparagraph{Case B}.  In this case, we assume that $\widehat{h} \geq
n^{-1/(2+\gamma)} (1 + \widehat{I} + I_0)^{1/2}$ and $\widehat{I} < 1
+ I_0$.  From equation~\eqref{Eqn-nicebound}, we have either
\begin{equation*}
\widehat{h}^2/4 + \lambda_n \widehat{I}^2/2 \leq O_\Pdist(n^{-1/2})
\widehat{h}^{1-\gamma/2}(1+I_0)^{1/2+ \gamma/4},
\quad \mbox{or} \quad
\widehat{h}^2/4 + \lambda_n \widehat{I}^2/2 \leq \lambda_n I_0^2/2.
\end{equation*}
These conditions imply, respectively, that
\begin{equation*}
\widehat{h} \leq (1+I_0)^{1/2}O_\Pdist(n^{-1/(\gamma+2)}), \qquad 
\mbox{and} \qquad
\widehat{I} \leq 1 + I_0
\end{equation*}
or
\begin{equation*}
\widehat{h}\leq O_\Pdist(\lambda_n^{1/2}I_0), \qquad \mbox{and} \qquad
\widehat{I} \leq O_\Pdist(I_0).
\end{equation*}
In either case, the proof is concluded by setting $\lambda_n^{-1} =
\order_\Pdist(n^{2/(\gamma+2)}(1+I_0))$. \\

\myparagraph{Case C}.  In this case, we assume that $\widehat{h} \leq
n^{-1/(2+\gamma)} (1 + \widehat{I} + I_0)^{1/2}$ and $\widehat{I} \geq
1 + I_0$.  From equation~\eqref{Eqn-nicebound}, we have
\begin{equation*}
\widehat{h}^2/4 + \lambda_n \widehat{I}^2/2 \leq
O_\Pdist(n^{-2/(2+\gamma)}) \widehat{I},
\end{equation*}
which implies that $\widehat{h} \leq
O_\Pdist(n^{-1/(2+\gamma)})\widehat{I}^{1/2}$ and $\widehat{I}\leq
\lambda_n^{-1} O_\Pdist(n^{-2/(2+\gamma)})$. Consequently, by setting
$\lambda_n^{-1} = \order_\Pdist(n^{2/(2+\gamma)})(1+I_0)$, we obtain
\begin{equation*}
\widehat{h} \leq O_\Pdist(\lambda_n^{1/2})(1+I_0), \qquad \mbox{and}
\qquad \widehat{I} \leq \order_\Pdist(1+I_0).
\end{equation*}

\vspace*{.02in}

\myparagraph{Case D}.  In this final case, we assume that $\widehat{h}
\leq n^{-1/(2+\gamma)} (1 + \widehat{I} + I_0)^{1/2}$ and $\widehat{I}
\leq 1 + I_0$, and the claim of Theorem~\ref{ThmConvE2}(a) follows. \\

\vspace*{.1in}

\newcommand{\Term}{\ensuremath{T}}

We now proceed to the proof of Theorem~\ref{ThmConvE2}(b).  Note that
part (a) and equation~\eqref{Ass-supnorm-penalize} imply that
$\|\gest\|_\infty = \order_\Pdist(1+I(\gtrue))$. Without loss of
generality, assume that $0< \eta_0 \leq \gtrue(x)$ and $\gest(x) \leq
\eta_1$ for all $x \in \Xspace$.  Then we have
\begin{eqnarray*}
\widehat{D}_K - D_K(\Pdist,\Qdist) & = & \int \log \gest
d\Pdist_n - \int \gest d\Qdist_n - (\int \log \gtrue d\Pdist - \int
\gtrue d\Qdist) \\ 
&\leq & \int \log \gest
d(\Pdist_n - \Pdist) - \int \gest d(\Qdist_n -\Qdist) =: T_1
\end{eqnarray*}
Also,
\begin{eqnarray*}
\widehat{D}_K - D_K(\Pdist,\Qdist)
& \geq & 
\int \log \gtrue d(\Pdist_n - \Pdist) - \int \gtrue d(\Qdist_n -\Qdist) =: T_2
\end{eqnarray*}

We have $\Term_2 = \order_\Pdist(n^{-1/2})$ by the central limit
theorem.   To bound $\Term_1$, we apply a modulus of continuity result
on the suprema of empirical processes with respect to function
$g$ and $\log g$, where $g$ is restricted to
smooth functions bounded from below (by $\eta_0$) and above (by
$\eta_1$).  As a result, the bracketing entropy for both function classes
$\Gclass$ and $\log \Gclass$ has the same order $O(1/\delta)^{\gamma}$
as given in~\eqref{Ass-entropy}. Apply Lemma 5.13 of~\cite{VdG-book}
(page 79), to obtain that for $\delta_n = n^{-1/(2+\gamma)}$, there
holds:
\begin{equation*}
\Term_1 = \order_\Pdist(n^{-1/2}\|\gest - \gtrue\|_{L_2(\Qdist)}^{1-
\gamma/2} \vee \delta_n^{2}) \; = \; \order_\Pdist(n^{-2/(2+\gamma)}),
\end{equation*}
thanks to part (a) of the theorem.
For $\gamma < 2$, $-2/(2+\gamma)< -1/2$. So the overall
rate is $n^{-1/2}$.

\section{Algorithmic implementation and simulation examples}
\label{Sec:algorithm}

In this section, we turn to the more practical issue of
implementation, focusing in particular on estimator E2.  When
$\Gclass$ has a kernel-based structure, we show how, via conversion to
its dual form, the computation of the estimator E2 reduces to the
solution of an $n$-dimensional convex program.  We illustrate the
performance of this practical estimator with a variety of simulations.

\subsection{Algorithms for kernel-based function classes}

We develop two versions of estimator E2: in the first, we assume that
$\Gclass$ is a reproducing kernel Hilbert space (RKHS), whereas in the
second, we assume that $\log \Gclass$ forms an RKHS.  In both cases,
we focus on the Hilbert space induced by a Gaussian kernel.  This
choice is appropriate as it is sufficiently rich, but also amenable to
efficient optimization procedures~\cite{Scholkopf02}.

We begin with some background on reproducing kernel Hilbert spaces;
see the books~\cite{Saitoh88,Scholkopf02} for further details.
Consider a positive definite function $\Ker$ mapping the Cartesian
product $\Xspace \times \Xspace$ to the non-negative reals.  By
Mercer's theorem, any such kernel function $\Ker$ can be associated
with a feature map $\Phi: \Xspace \rightarrow \mathcal{H}$, where
$\mathcal{H}$ is a Hilbert space with inner product $\inprod{\cdot}{\cdot}$.
Moreover, for all $x,x'\in
\Xspace$, the inner product in this Hilbert space is linked to the
kernel via the relation $\Ker(x,x') = \inprod{\Phi(x)}{\Phi(x')}$. As
a reproducing kernel Hilbert space, any function $g \in \mathcal{H}$
can be expressed as an inner product $g(x) = \inprod{\w}{\Phi(x)}$,
where $\|g\|_\Hspace = \|\w\|_{\Hspace}$. The kernel used in our
simulations is the Gaussian kernel:
\begin{eqnarray}
\label{EqnGaussKernel}
K(x,y) & \defn & \exp \left \{ -\|x - y\|^2/\sigma \right \},
\end{eqnarray}
where $\| \cdot \|$ is the Euclidean metric in $\real^\dims$, and $\sigma >
0$ is a parameter for the function class.

\subsubsection{Imposing RKHS structure of $\Gclass$}

Suppose that the function class $\Gclass$ in estimator E2 is the
Gaussian RKHS space $\Hspace$, and let the complexity measure be the
Hilbert space norm $I(g) = \|g\|_\Hspace$.  With these choices,
equation~\eqref{Eqn-penalize} becomes:
\begin{equation}
\label{Eqn-program-G}
f^* \; = \; \min_{\w} J \defn \min_{\w} \frac{1}{n}\sum_{i=1}^{n}
\inprod{\w}{\Phi(x_i)} - \frac{1}{n}\sum_{j=1}^{n}\log
\inprod{\w}{\Phi(y_j)} + \frac{\lambda_n}{2}\|w\|_\Hspace^2,
\end{equation}
where the samples $\{x_i\}$ and $\{y_j\}$ are i.i.d.  from $\Qdist$
and $\Pdist$, respectively. The $\log$ function is extended to take
the value $-\infty$ for negative arguments.

\begin{lemma}
\label{Lem-dual-G} The primal value $f^*$ has the equivalent
dual form:
\begin{multline*}
\label{Eqn-dual-G}
-\min_{\alpha > 0} \; \left \{-1 -\frac{1}{n} \sum_{j=1}^{n} \log
n\alpha_j + \frac{1}{2\lambda_n} \sum_{i,j}\alpha_i \alpha_j
K(y_i,y_j) + \frac{1}{2\lambda_n n^2}\sum_{i,j}K(x_i,x_j) -
\frac{1}{\lambda_n n}\sum_{i,j}\alpha_j K(x_i,y_j) \right \}.
\end{multline*}
Moreover, the optimal dual solution $\alphaest$  is linked to
the optimal primal solution $\west$ via
the relation
\begin{eqnarray}
\label{EqnPrimalDual}
\west & = & \frac{1}{\lambda_n}(\sum_{j=1}^{n}\alphaest_j \Phi(y_j) -
\frac{1}{n}\sum_{i=1}^{n}\Phi(x_i)).
\end{eqnarray}
\end{lemma}
\begin{proof}
Let $\psi_i(\w) \defn \frac{1}{n}\inprod{\w}{\Phi(x_i)}$,
$\varphi_j(\w) \defn -\frac{1}{n}\log \inprod{\w}{\Phi(y_j)}$, and
$\Omega(\w) = \frac{\lambda_n}{2}\|\w\|^2_\Hspace$.  We have
\begin{eqnarray*}
f^* & = & -\max_{\w}(\inprod{0}{\w} - J(\w)) = -J^*(0) \\ & = &
-\min_{u_i,v_j} \sum_{i=1}^{n}\psi_i^*(u_i) +
\sum_{j=1}^{n}\varphi_j^*(v_j) + \Omega^*(-\sum_{i=1}^{n} u_i -
\sum_{j=1}^{n} v_j),
\end{eqnarray*}
where the last line is due to the inf-convolution
theorem~\cite{Rockafellar}.  Simple calculations yield:
\begin{eqnarray*}
\varphi_j^*(v) & = & \begin{cases} -\frac{1}{n} - \frac{1}{n}\log n\alpha_j
& \mbox{if} \; v=-\alpha_j \Phi(y_j) \\
 +\infty
& \mbox{otherwise}   \end{cases} \\
\psi_i^*(u) & = & 0 \;\mbox{if}\; 
u = \frac{1}{n} \Phi(x_i) \;\mbox{and}\; +\infty \;\mbox{otherwise} \\
\Omega^*(v) & = & \frac{1}{2\lambda_n}\|v\|_\Hspace^2.
\end{eqnarray*}
Thus, we conclude that $f^* = -\min_{\alpha_i}\sum_{j=1}^{n}
(-\frac{1}{n} - \frac{1}{n}\log n\alpha_j) + \frac{1}{2\lambda_n}
\|\sum_{j=1}^{n}\alpha_j\Phi(y_j) - \frac{1}{n}\sum_{i=1}^{n}
\Phi(x_i)\|_\Hspace^2$, from which the claim follows.
The primal-dual relation~\eqref{EqnPrimalDual} also follows from
these calculations.
\end{proof}

For an RKHS based on a Gaussian kernel, the entropy
condition~\eqref{Ass-entropy} holds for any $\gamma > 0$
(cf.\ Zhou~\cite{Zhou}).  Furthermore,
condition~\eqref{Ass-supnorm-penalize} holds since, via the Cauchy-Schwarz
inequality, we have
\begin{equation*}
|g(x)| = |\inprod{\w}{\Phi(x)}|\leq \|w\|_\Hspace
\|\Phi(x)\|_\Hspace \leq I(g) \sqrt{K(x,x)} \leq I(g).
\end{equation*}
Thus, by Theorem~\ref{ThmConvE2}(a), we have $\|\west\|_\Hspace =
\|\gest\|_\Hspace = \order_\Pdist(\|\gtrue\|_\Hspace)$, so the penalty
term $\lambda_n \|\west\|^2$ vanishes at the same rate as $\lambda_n$.
Thus, we obtain the following estimator for the KL
divergence:
\begin{equation}
\label{EqnGclassEst}
\widehat{D}_K = 1 + \sum_{j=1}^{n} (-\frac{1}{n} - \frac{1}{n}\log
n\alphaest_j) = \sum_{j=1}^{n} -\frac{1}{n}\log n\alphaest_j.
\end{equation}


\subsubsection{Imposing RKHS structure on $\log \Gclass$}  

An alternative starting point is to posit that the function class $\log
\Gclass$ has an RKHS structure.  In this case, we consider functions
of the form $g(x) = \exp \inprod{\w}{\Phi(x)}$, and use the
complexity measure $I(g) = \|\log g\|_\Hspace = \|w\|_\Hspace$.
Unfortunately, Theorem~\ref{ThmConvE2} does not apply directly because
condition~\eqref{Ass-supnorm-penalize} no longer holds, but this
choice nonetheless seems reasonable and worth investigating from an
empirical viewpoint.

A derivation similar to the previous case yields the following convex
program:
\begin{eqnarray*}
\label{Eqn-program}
\min_{\w} J & \defn & \min_{\w} \;
\frac{1}{n}\sum_{i=1}^{n}e^{\inprod{\w}{\Phi(x_i)}} -
\frac{1}{n}\sum_{j=1}^{n}\inprod{\w}{\Phi(y_j)} +
\frac{\lambda_n}{2}\|w\|_\Hspace^2 \\ & = & -\min_{\alpha > 0} \;
\sum_{i=1}^{n}\alpha_i \log(n\alpha_i) -\alpha_i +
\frac{1}{2\lambda_n} \|\sum_{i=1}^{n}\alpha_i\Phi(x_i) -
\frac{1}{n}\sum_{j=1}^{n} \Phi(y_j)\|_\Hspace^2.  \comment{ & = &
-\min_{\alpha \geq 0} \; \sum_{i=1}^{n} \biggr (\alpha_i\log \alpha_i
+ \alpha_i \log\frac{n}{e} \biggr ) + \frac{1}{2\lambda_n}\sum_{1\leq
i,j\leq n}\alpha_i K(x_i,x_j)\alpha_j \\ & - & \frac{1}{\lambda_n
n}\sum_{1\leq i,j\leq n}\alpha_i K(x_i,y_j) + \frac{1}{2\lambda_n
n^2}\sum_{1\leq i,j\leq n}K(y_i,y_j).}
\end{eqnarray*}
Letting $\alphaest$ be the solution of the above convex program, the
KL divergence can be estimated by:
\begin{equation}
\label{EqnLogGclassEst}
\widehat{D}_K = 1 + \sum_{i=1}^{n}\alphaest_i\log \alphaest_i +
\alphaest_i\log\frac{n}{e}.
\end{equation}

\subsection{Simulation results}
\label{Sec:simulation}

In this section, we describe the results of various simulations that
demonstrate the practical viability of the
estimators~\eqref{EqnGclassEst} and~\eqref{EqnLogGclassEst}, as well
as their convergence behavior.  We experimented with our estimators
using various choices of $\Pdist$ and $\Qdist$, including Gaussian,
beta, mixture of Gaussians, and multivariate Gaussian distributions.
Here we report results in terms of KL estimation error.  For each of
the eight estimation problems described here, we experiment with
increasing sample sizes (the sample size, $n$, ranges from $100$ to
$10^4$ or more). Error bars are obtained by replicating each set-up
250 times.

For all simulations, we report our estimator's performance using the
simple fixed rate $\lambda_n \sim 1/n$, noting that this may be a
suboptimal rate. We set the kernel width to be relatively small
($\sigma = .1$) for one-dimensional data, and choose larger $\sigma$
for higher dimensions. We use M1 and M2 to denote the
estimators~\eqref{EqnGclassEst} and~\eqref{EqnLogGclassEst},
respectively.  We compare these methods to algorithm $A$ in Wang et
al~\cite{Wang05}, which was shown empirically to be one of the best
methods in the literature.  This method, to be denoted by WKV, is
based on data-dependent partitioning of the covariate
space. Naturally, the performance of WKV is critically dependent on
the amount $s$ of data allocated to each partition; here we report
results with $s \sim n^{\gamma}$, where $\gamma = 1/3,1/2,2/3$.

\begin{figure}
{\centering
\begin{tabular}{cc}
\includegraphics[width=.45\textwidth,height=.40\textwidth]
{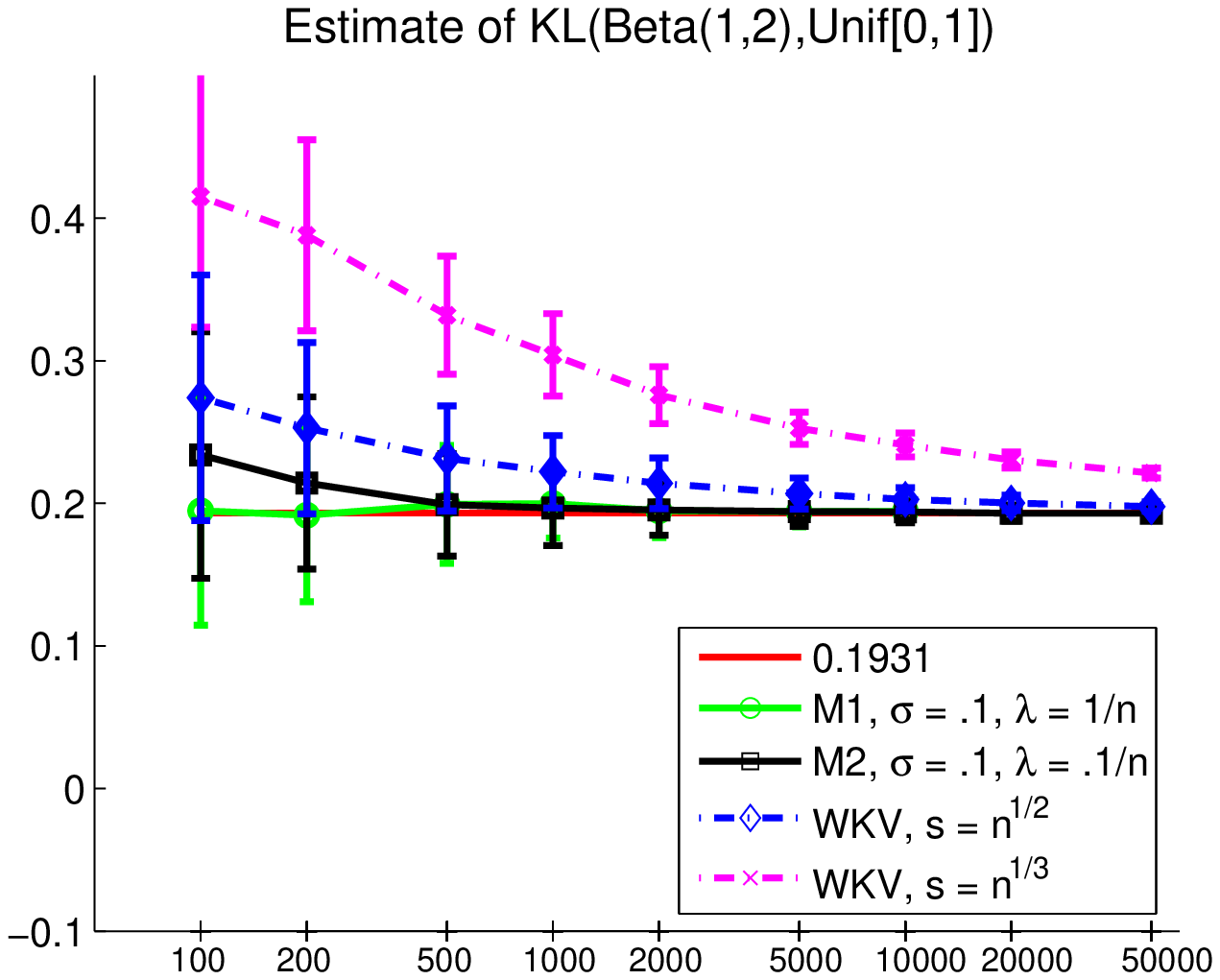} &
\includegraphics[width=.45\textwidth,height=.40\textwidth]
{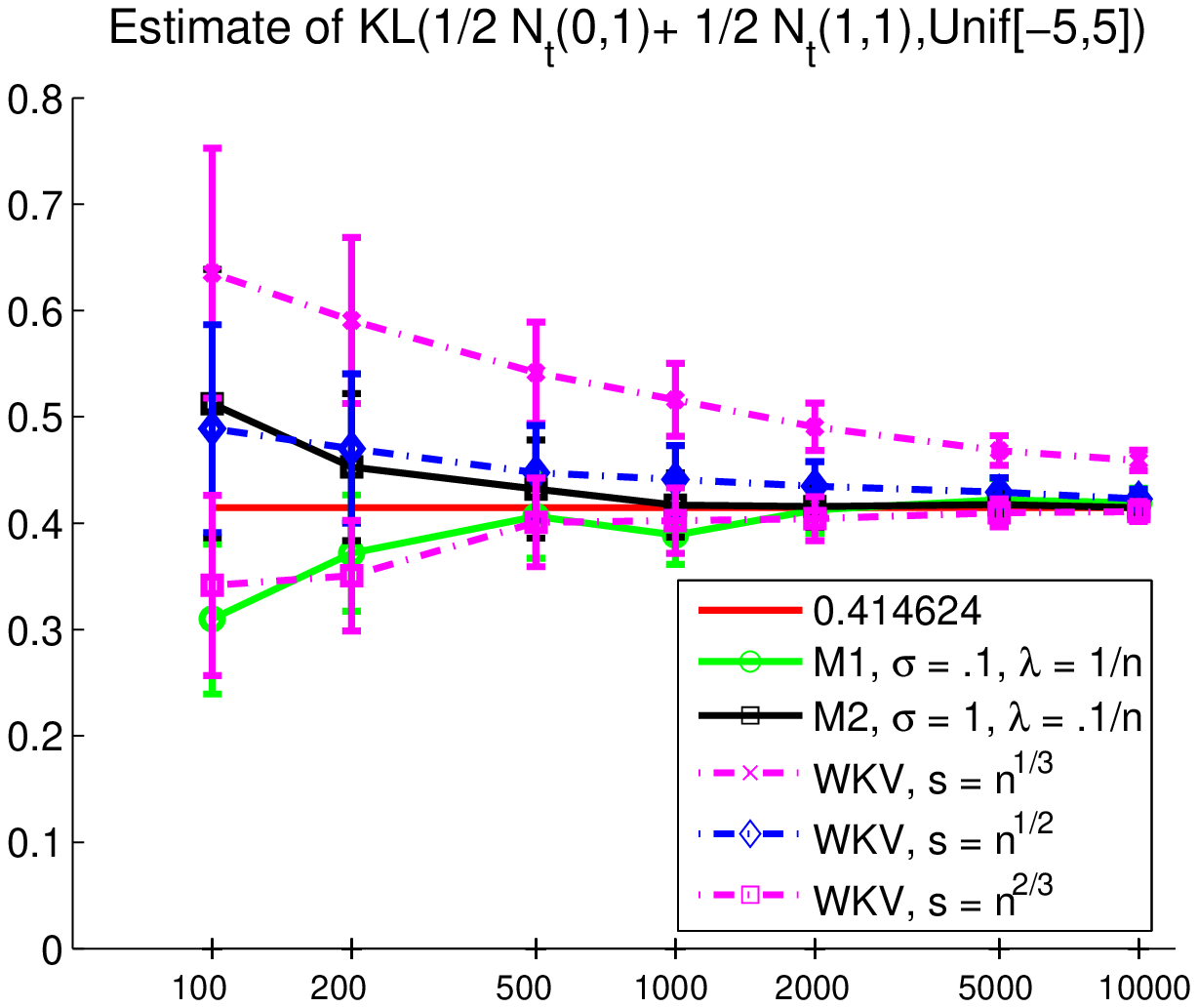}\\
\includegraphics[width=.45\textwidth,height=.40\textwidth]
{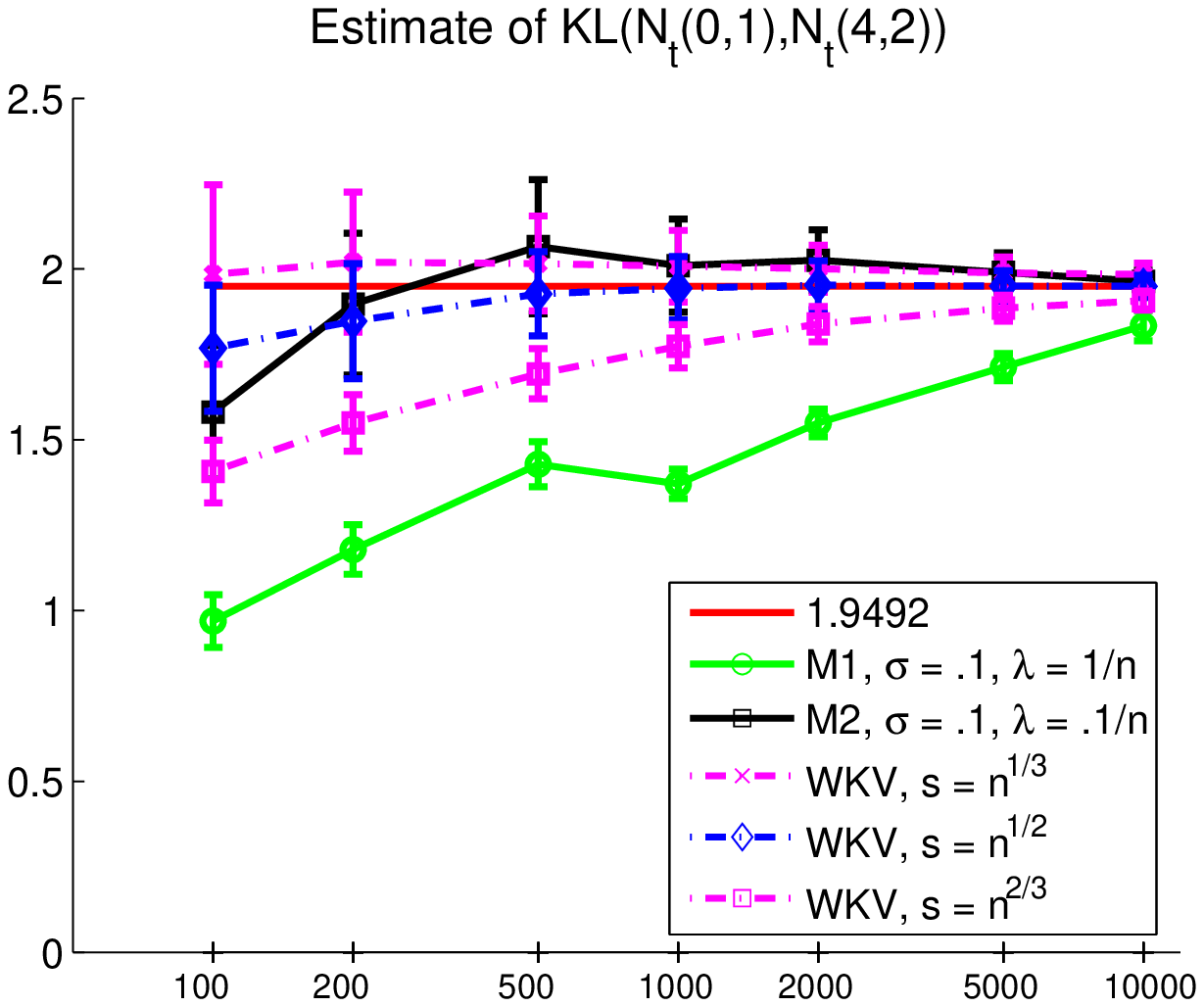} &
\includegraphics[width=.45\textwidth,height=.40\textwidth]
{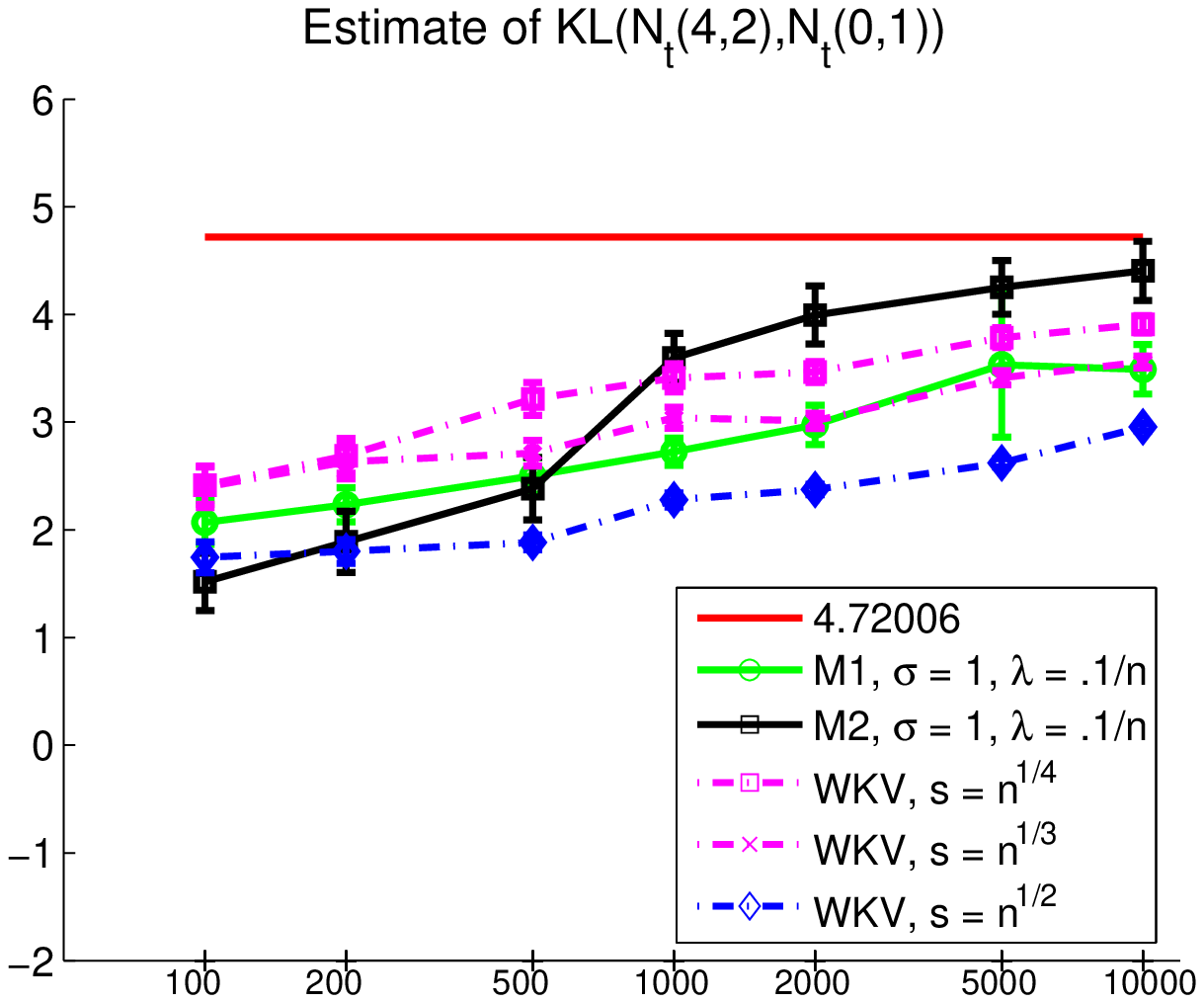} \\
\end{tabular}
\caption{Results of estimating KL divergences for various choices
of probability distributions. In all plots, the X-axis is the number 
of data points plotted on a log scale, and the Y-axis is the estimated value.
The error bar is obtained by replicating the simulation 250 times.
$N_t(a,I_k)$ denotes a truncated normal distribution of
$k$ dimensions with mean $(a,\ldots,a)$ and identity covariance
matrix.}
}
\label{fig:univariate}
\end{figure}

\begin{figure}
{\centering
\begin{tabular}{cc}
\includegraphics[width=.45\textwidth,height=.40\textwidth]
{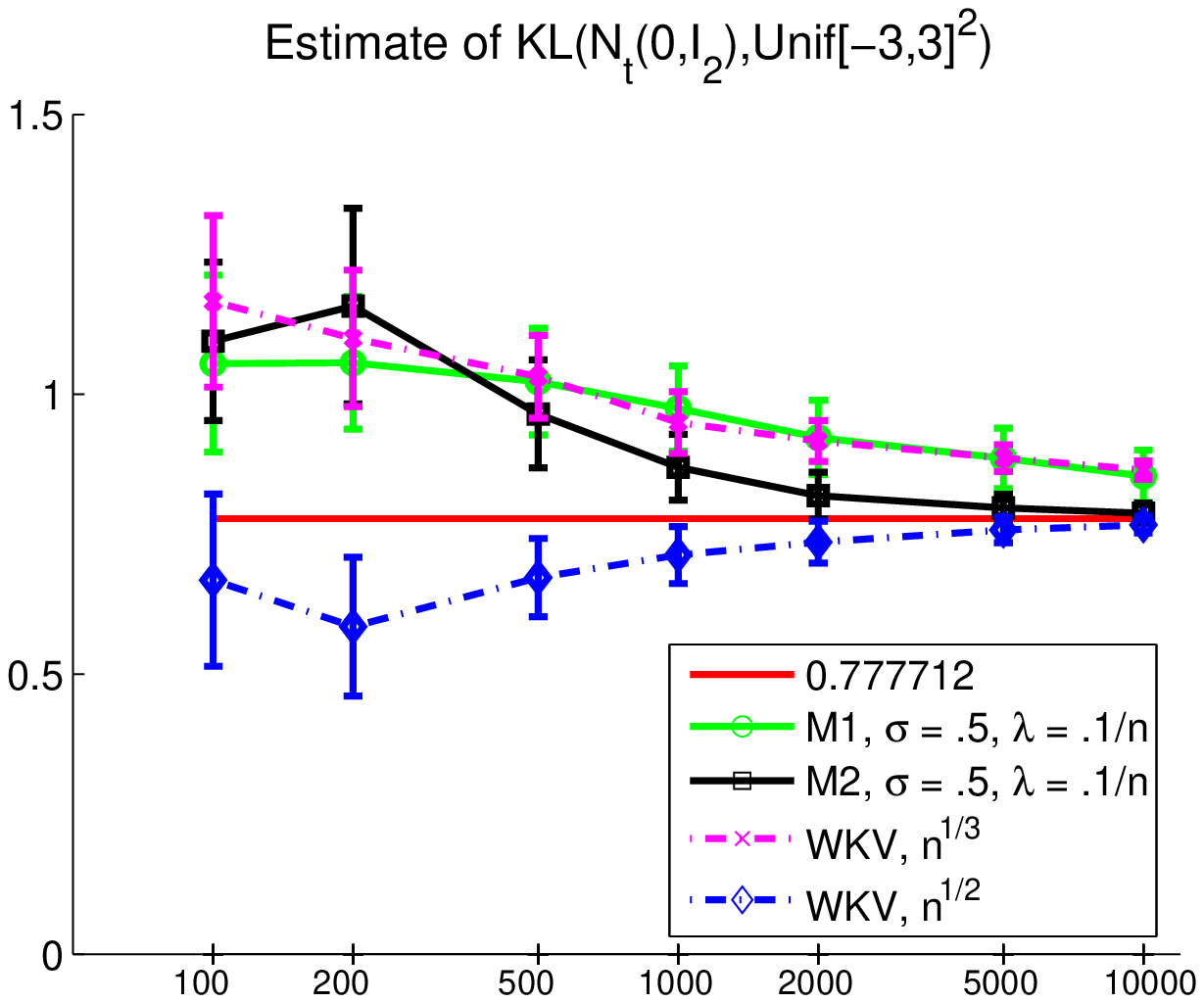} &
\includegraphics[width=.45\textwidth,height=.40\textwidth]
{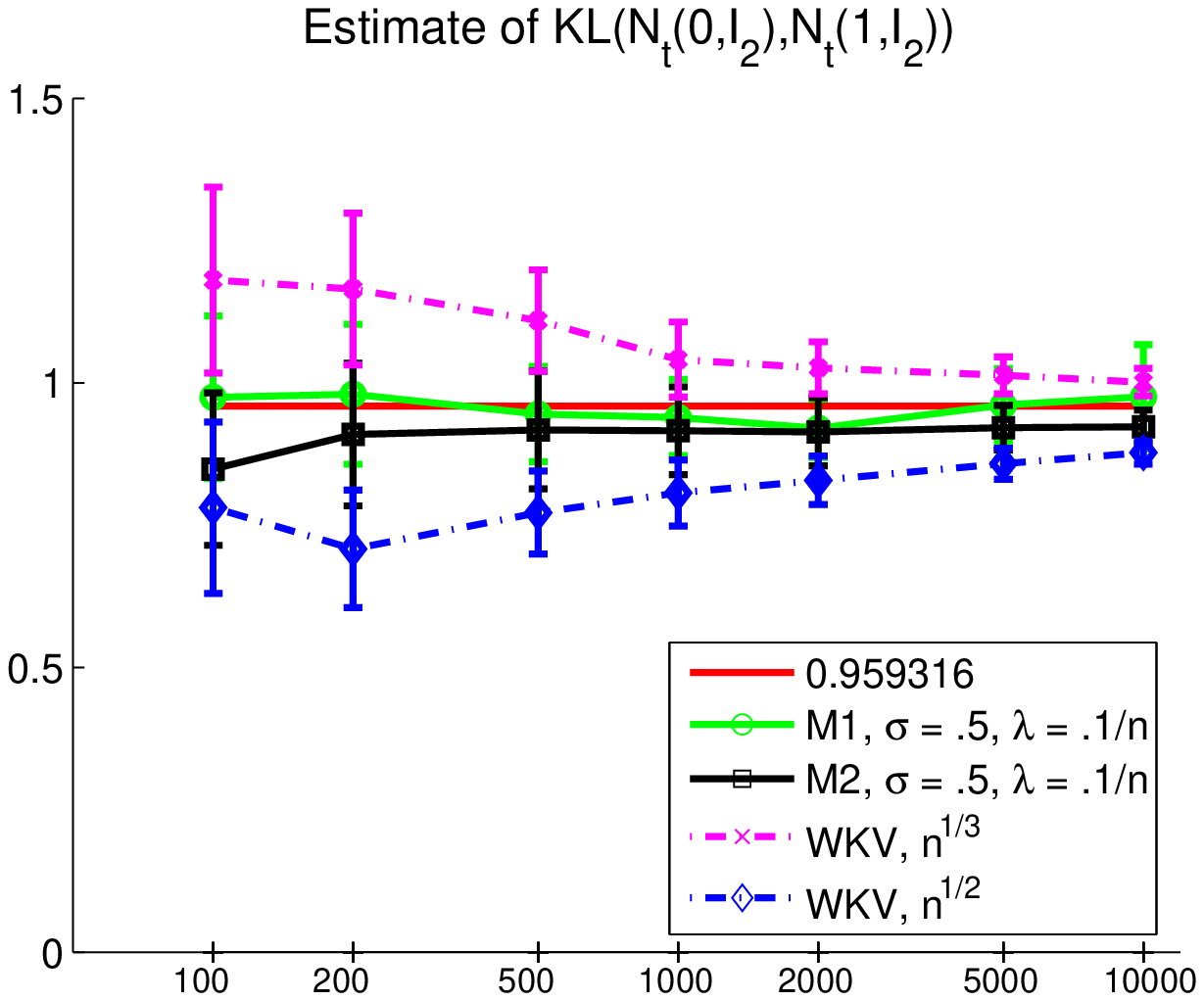} \\
\includegraphics[width=.45\textwidth,height=.40\textwidth]
{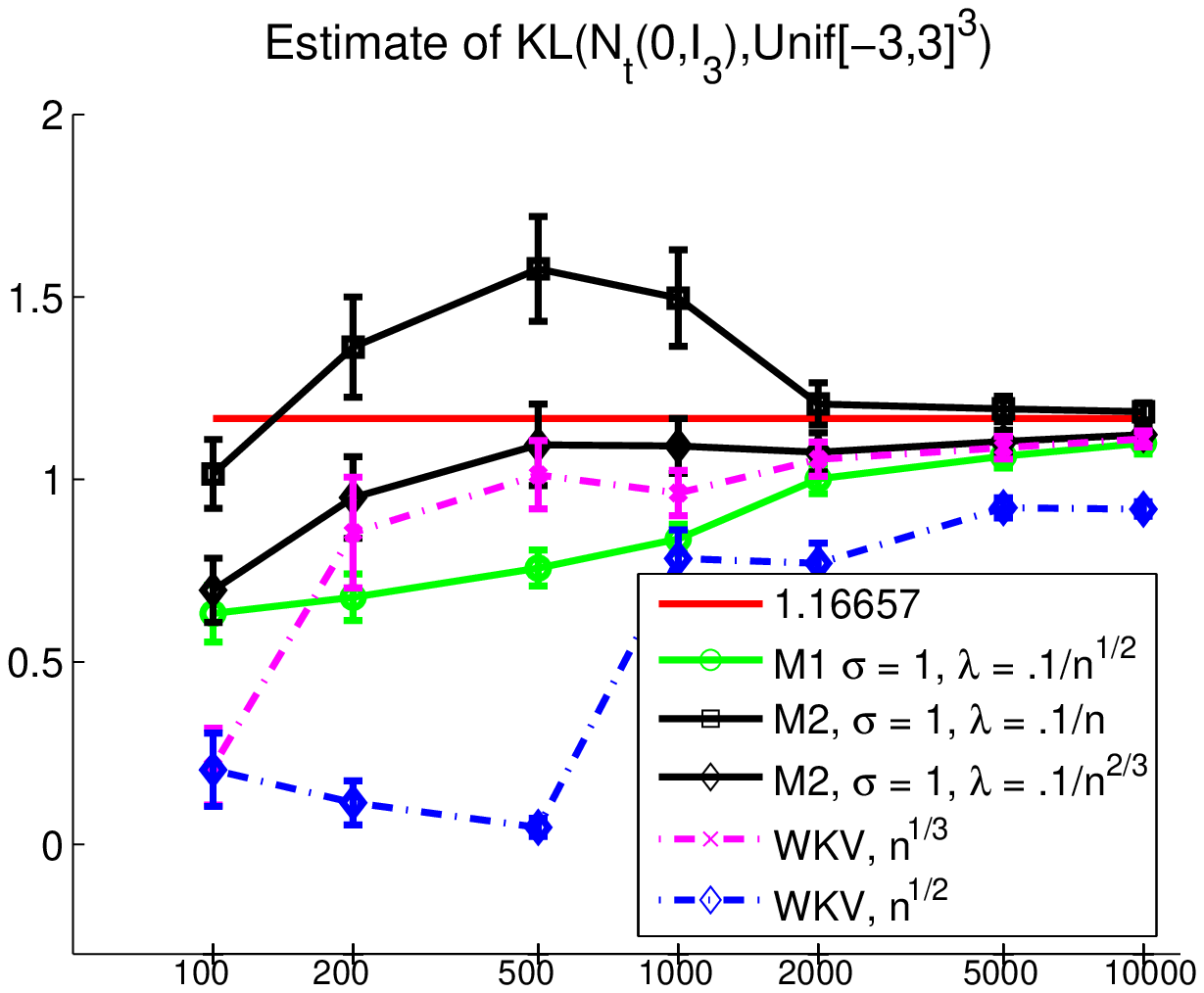} &
\includegraphics[width=.45\textwidth,height=.40\textwidth]
{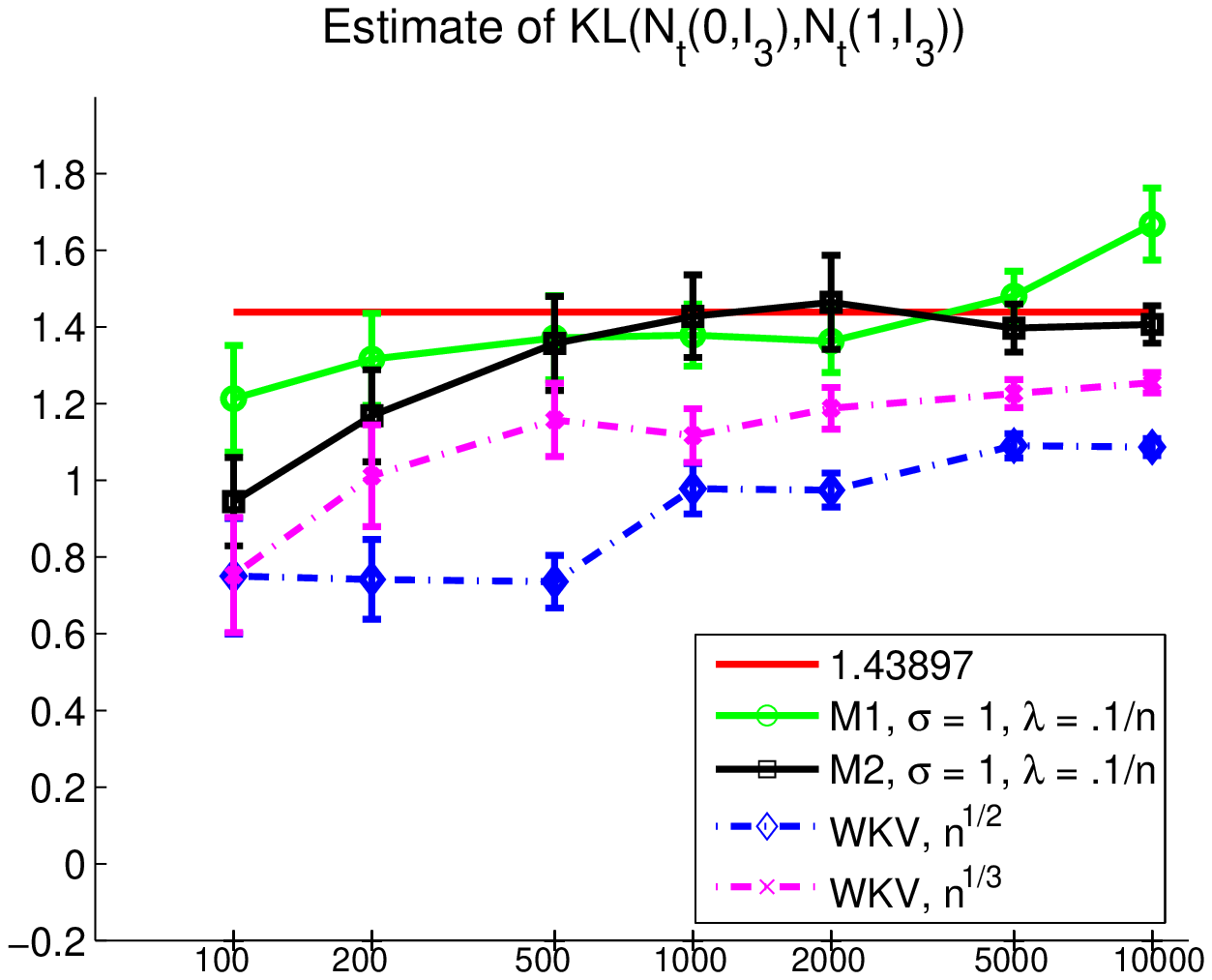} 
\end{tabular}
\caption{Results of estimating KL divergences for various choices
of probability distributions. In all plots, the X-axis is the number 
of data points plotted on a log scale, and the Y-axis is the estimated value.
The error bar is obtained by replicating the simulation 250 times.
$N_t(a,I_k)$ denotes a truncated normal distribution of
$k$ dimensions with mean $(a,\ldots,a)$ and identity covariance
matrix.}
}
\label{fig:multivariate}
\end{figure}

The four plots in Figure~\ref{fig:univariate} present results 
with univariate distributions.  We see that the estimator M2
generally exhibits the best convergence rate among the 
estimators considered.  The WKV estimator performs somewhat 
less well, and shows sensitivity to the choice of partition
size $s$, with the ranking of the different WKV estimators 
changing over the experiments.  The performance of estimator 
M1 is comparable to that of the WKV estimator, although clearly 
better in the first plot.  In Figure~\ref{fig:multivariate} 
we present the results with two- and three-dimensional data. 
Again, estimator M2 has the best convergence rates in all 
examples. The M1 estimator does not converge in the last 
example, suggesting that the underlying function class exhibits 
very strong bias.  In these examples, the WKV estimator again
shows sensitivity to the choice of partition size; moreover,
its performance is noticeably degraded in the case of 
three-dimensional data (the lower two plots).

It is worth noting that as one increases the number of dimensions, 
histogram-based methods such as WKV become increasingly difficult 
to implement, whereas increased dimension has only a mild effect 
on the complexity of implementation of our method.

\section{Some extensions}
\label{Sec:general}

In this section, we discuss some extensions and related estimators,
all based on the same basic variational principle.

\subsection{Estimation of likelihood ratio functions}

Suppose that we are primarily interested in estimating the likelihood
ratio function $\ptrue/\qtrue$, as opposed to the Kullback-Leibler
divergence.  In this case, we may consider any divergence functional
$D_\phi(\ptrue,\qtrue)$, where $\phi$ is a convex function on
$\real_{+}$, possibly different than the logarithm leading to KL
divergence.  Again applying Lemma~\ref{Lem-variational}, choosing a
different divergence leads to the following alternative estimator of
the likelihood ratio:
\begin{eqnarray}
\label{Eqn-LR}
\fest & \defn & \argmax_{f\in \Fclass} \int f\; d\Qdist_n - \int
\phi^*(f) \;d\Pdist_n \\
\label{Eqn-Dphi}
\widehat{D}_\phi & \defn & \int \fest\; d\Qdist_n - \int \phi^*(\fest)
\;d\Pdist_n.
\end{eqnarray}
The quantity $\fest$ is an estimate of the quantity $\ftrue = \partial
\phi(\qtrue/\ptrue)$, whereas $\widehat{D}_\phi$ is an estimate of the
divergence $D_\phi(\ptrue, \qtrue)$ (of secondary interest for the
moment).

We make the following observations:
\begin{itemize}
\item If $\phi$ is a differentiable and strictly convex function,
i.e., $\phi'' > 0$, then the likelihood ratio function $\ptrue/\qtrue$
can be recovered by applying the inverse of $\phi'$ to $\fest$. Thus,
we obtain a family of estimation methods for the likelihood ratio
function by simply ranging over choices of $\phi$.

\item If (on the other hand) the function $\phi$ is chosen to be
non-differentiable, we cannot directly invert the mapping $\partial
\phi$, but we can nonetheless obtain estimators for other interesting
objects. For instance, suppose that $\phi$ has the piecewise-linear
form
\begin{eqnarray*}
\phi(u) & \defn &  \begin{cases}|u-1| \;\mbox{if}\; u \geq 0 \\ +\infty \;
\mbox{otherwise},
\end{cases},
\end{eqnarray*}
so that $D_\phi$ is the variational distance. Noting that $\partial
\phi(u) = \sign(u-1)$ for any $u > 0$,  we see that the quantity
$\fest$ in equation~\eqref{Eqn-LR} provides an estimate of the
thresholded likelihood ratio.\footnote{
In fact, there is strong relationship between variational distance and
a threshold function of the likelihood ratio.  Note that the conjugate
dual for $\phi$ has the form:
\[\phi^*(v) = \begin{cases}-1 \;\mbox{if}\; v \geq -1 \\ 
v \; \mbox{if}\; v \in [-1,1] \\ +\infty \;\mbox{otherwise,}
\end{cases}\]
which is related to a hinge loss in the literature of binary
classification in machine learning. Indeed, a binary classification
problem can be viewed as estimating the threshold function of the likelihood
ratio. See~\cite{NguyenEtAl-AOS} for a discussion of divergences and
surrogate losses from this viewpoint.}
\end{itemize}

\subsection{Extensions to different $\phi$}

Let us assume that $\phi$ is chosen to be differentiable and strictly
convex, so that we can estimate the likelihood ratio $\gtrue$ by
applying $(\phi')^{-1}$ to $\fest$.  Since there are many such $\phi$,
it is natural to ask how the properties of $\phi$ affect the quality
of the estimate of $\gtrue$.  The analysis provided in the preceding
sections can be adapted to other choices of $\phi$, as we describe
here.

In order to describe these extensions, we first define a distance
between $\ftrue$ and $f$:
\begin{eqnarray}
\label{Eqn-distance-phi}
d_\phi(\ftrue,f) & \defn & D_\phi(\Pdist,\Qdist) - \int f d\Qdist -
 \phi^*(f)d\Pdist
\end{eqnarray}
Note that this distance is simply the generalization of the quantity
$d(\gtrue,g)$ previously defined in Lem~\ref{Prop-dist-KL}).  For
future reference, we note the equivalence
\begin{eqnarray*}
d_\phi(\ftrue,f) & = &\int (\phi^*(f) - \phi^*(\ftrue)) d\Pdist - (f -
 \ftrue)d\Qdist \\ & = & \int \left(\phi^*(f) - \phi^*(\ftrue) -
 \frac{\partial \phi^{*}}{\partial f} (\ftrue) (f-\ftrue) \right) \;d\Pdist,
\end{eqnarray*}
where the final line uses the facts that $\frac{\partial
\phi^{*}}{\partial f}(\ftrue) = \qtrue/\ptrue$ and
$\phi'(\qtrue/\ptrue) = \ftrue$.  This expression shows that $d_\phi$
is the \emph{Bregman divergence} defined by the convex function
$\phi^*$.

Recall that the key ingredient in our earlier analysis was the
relation between the empirical processes defined by
equation~\eqref{Eqn-supreme} and the ``distance'' $d(\gtrue,g)$ (see
Lemma~\ref{Prop-dist-KL}).  Similarly, the key technical ingredient in
the extension to general $\phi$ involves relating
the quantity
\begin{eqnarray*}
v_n^\phi(\Fclass) & = &
\sup_{f\in \Fclass}
\biggr |\int (\phi^*(f) - \phi^*(\ftrue)) d(\Pdist_n - \Pdist) -
\int (f - \ftrue) d(\Qdist_n - \Qdist) \biggr |
\end{eqnarray*}
to the distance $d_\phi(\ftrue,f)$ defined in
equation~\eqref{Eqn-distance-phi}.  In particular, we can state the
following analog of Lemma~\ref{Prop-dist-KL} and Lemma~\ref{Lem-hellinger}:
\begin{lemma}
\label{LemGeneral}
Let $\fest$ be the estimate of $\ftrue$ obtained by solving the
problem~\eqref{Eqn-LR}.  Then
\begin{equation}
d_\phi(\ftrue,\fest) \; \leq \; v_n^\phi(\Fclass).
\end{equation}
\end{lemma}

Under suitable technical conditions, we have $v_n^\phi(\Fclass)
\convas 0$, so that Lemma~\ref{LemGeneral} implies that $\fest$ is a
consistent estimator for $\ftrue$ in the sense of $d_\phi$. This lemma
also provides the technical means to derive convergence rates in the
same manner as in the previous sections.  Note
that $d_\phi(\ftrue,f)$ is usually not a proper metric. To apply
standard results from empirical process theory, the trick is that one
can replace $d_\phi$ by a lower bound which is a proper metric (such
as $L_2$ or Hellinger metric). In the case of KL divergence, we have
seen that this lower bound is the Hellinger distance (via
Lemma~\ref{Prop-dist-KL}(i)).

Let us illustrate this idea by stating a result about likelihood ratio
estimation in terms of the \emph{$\chi$-square divergence},
which is defined by
\begin{eqnarray}
\label{EqnDefnChiSquare}
D_\chi(\Pdist,\Qdist) & \defn & \int
\ptrue^2/\qtrue d\mu.
\end{eqnarray}
Note that this divergence is an $f$-divergence with $\phi(u) = 1/u$; 
a short calculation (see Appendix~\ref{AppCalc}) shows that the associated
``distance'' $d_\chi$ is given by $d_\chi(g, \gtrue) = \int(g -
\gtrue)^2 d\Qdist$, which is simply the $L_2(\Qdist)$ metric.
With this set-up, the estimator now has the following ``least square''
form:
\[
\gest \defn \argmax_{g\in \Gclass} \int -g^2\; d\Qdist_n + \int
2g \;d\Pdist_n 
\]
The following theorem is an analog of 
Theorem~\ref{ThmConvE1} (with an almost identical proof):
\begin{theorem}
\label{ThmExt}
Assume that for some constant $0 < \gamma < 2$,
\begin{eqnarray}
\label{Eqn-bounded-entropy}
\entrobrack{\delta}{\Gclass}{L_2(\Qdist)} & \leq & A_\Gclass \delta^{-\gamma},
\end{eqnarray}
and moreover that condition~\eqref{Eqn-Upbound} holds.  Then the
estimator $\gest$ obtained from the $\chi$-square divergence is
consistent with rate $d_\chi(\gtrue,\gest) = O_\Pdist(n^{-2/(\gamma +
2)})$.
\end{theorem}

\vspace*{.05in}

\myparagraph{Remark:} Comparing with Theorem~\ref{ThmConvE1}, we see
that the conditions of Theorem~\ref{ThmExt} are weaker.  Indeed, the
$L_2(\Qdist)$ metric is dominated by the Hellinger metric, so that
imposing bounds on $L_2(\Qdist)$-metric and its induced entropy are
milder conditions.


\subsection{Estimation of the divergence functional $D_\phi$}


Suppose that we are primarily interested in estimating the divergence
functional $D_\phi$, given that we have already obtained an optimal
estimator $\gest$ of the likelihood ratio function $\gtrue =
\ptrue/\gtrue$ (such as the one defined by~\eqref{Eqn-g}
or~\eqref{Eqn-penalize}, or more generally~\eqref{Eqn-LR}). 
We have demonstrated that $D_\phi$ can be estimated
by~\eqref{Eqn-KLest} and \eqref{Eqn-penalize-KLest}, or more generally
by~\eqref{Eqn-Dphi}.  Note that $D_\phi$ can be viewed as an integral
of the likelihood ratio under the distribution $\Qdist$. Indeed, we
can write
\[D_\phi(\Pdist,\Qdist) = \int (\ptrue/\qtrue)\phi(\qtrue/\ptrue) \; d\Qdist
= \int \gtrue \phi(1/\gtrue) \; d\Qdist.\] 
Although $D_\phi$ is an integral functional of $\gtrue =
\ptrue/\qtrue$, an interesting feature here is that the integration is
with respect to \emph{unknown} $\Qdist$. In this section, we show that
estimators such as~\eqref{Eqn-KLest} and~\eqref{Eqn-penalize-KLest} for
the KL divergence can be viewed as a first-order Taylor expansion of
the integral functional around the estimate $\gest$ of the likelihood 
ratio. This discussion is motivated by a line of work on the
estimation of integral functional of a single density function
(cf.~\cite{Joe-1989,Birge-Massart}), and also leads to an open
question.

Suppose that $\phi:\real_+ \rightarrow \real$ is a convex function 
differentiable up to third order, $\Gclass$ is a smooth function
class bounded from both above and below as in Example 1 (with smoothness
parameter $\alpha$).  Suppose that $\gest$ is an estimator of $\gtrue$
(such as the one defined by~\eqref{Eqn-LR}), and that $\|\gest -
\gtrue\|_{L_2(\Qdist)} = O_P(n^{-\alpha/(2\alpha + \dims)})$.  Using a
Taylor expansion around $\gest$, we obtain:
\begin{eqnarray*}
g\phi(1/g) & = &\gest \phi(1/\gest) + (g-\gest)(\phi(1/\gest) -
\phi'(1/\gest)/\gest) + (g-\gest)^2 \phi''(1/\gest)/\gest^3 + \\ &&
O((g-\gest)^3) \\ & = & \phi'(1/\gest) + \phi''(1/\gest)/\gest +
g(\phi(1/\gest) - \phi'(1/\gest)/\gest - 2\phi''(1/\gest)/\gest^2) +
\\ && g^2 \phi''(1/\gest)/\gest^3 + O((g-\gest)^3).
\end{eqnarray*}
We arrive at
\begin{eqnarray*}
D_\phi(\Pdist,\Qdist) & = &\int g\phi(1/g) d \Qdist \\ 
& = & \int \left( \phi'(1/\gest) + \phi''(1/\gest)/\gest \right) \; d \Qdist \\ 
&& + \int \left(\phi(1/\gest) - \phi'(1/\gest)/\gest - 
2\phi''(1/\gest)/\gest^2\right) \; d\Pdist \\ 
&& + \int \ptrue^2/\qtrue \; \phi''(1/\gest)/\gest^3 \; d\mu +
O(\|\gtrue-\gest\|_3^3).
\end{eqnarray*}
In the above expression, the first two integrals can be estimated from
two $n$-samples of empirical data drawn from $\Pdist$ and $\Qdist$.
Because of the boundedness assumption, these estimations have at most
$O_P(n^{-1/2})$ error. The error of our Taylor approximation is
$O(\|\gtrue-\gest\|_3^3) = O_P(n^{-3\alpha/(2\alpha + \dims)})$.  This
rate is less than $O(n^{-1/2})$ for $\alpha \geq \dims/4$.  Thus when
$\alpha \geq \dims/4$, the optimal rate of convergence for estimating
$D_\phi$ hinges on the convergence rate for estimating the integral of the form
$\int \ptrue^2/\qtrue \psi \; d\mu$. This is interesting because
we have reduced the problem of estimating any $f$-divergence to
a particular integral of two densities $\int \ptrue^2/\qtrue \psi \; d\mu$,
where $\psi$ is a known function.

\comment{ \myparagraph{Remark.}  Our method described in the previous
section yields the optimal rate of $O(n^{-1/2})$ for estimating $\int
\ptrue^2/\qtrue \psi \; d\mu$, given that $\alpha > d/2$.  We do not
know at the moment if this rate is attainable when $\alpha \in (d/4,
d/2]$, as achieved by~\citep{Bickel-Ritov} and~\citep{Birge-Massart}
for estimating integral functionals of \emph{one} density.  For
$\alpha < d/4$, the best minimax bound for estimating $D_\phi$ must be
no greater than the minimax lower bound for estimating integral
functional of one density, which was shown to be $n^{-4\alpha/(4\alpha
+ \dims)}$ through the work of~\cite{Bickel-Ritov},
~\cite{Birge-Massart} and~\cite{Kerkyacharian-Picard}. We do not know
at the moment whether this minimax lower bound is attainable in the
estimation of $D_\phi$.  }

Let us return to the case of KL divergence, i.e., $\phi(u) = -\log u$. If we 
use Taylor approximation up to first order, the estimator has
the following form:
\begin{eqnarray*}
\widehat{D}_\phi & = & \int (\phi(1/\gest) - \phi'(1/\gest)/\gest) \;
d\Pdist_n + \int \phi'(1/\gest)\; d\Qdist_n \\ & = & \int \log \gest +
1 d\Pdist_n - \gest d\Qdist_n,
\end{eqnarray*}
which has exactly the same form as our original
estimator~\eqref{Eqn-KLest}, except that here $\gest$ can be any
(optimal) estimator of the likelihood ratio. Note that we
have shown the optimal convergence rate of $n^{-1/2}$ for 
the KL divergence estimator, given $\alpha > \dims/2$ (so
that $\gamma = \dims/\alpha < 2$). Questions regarding the
estimator and its analysis for the case $\alpha \leq \dims/2$ remain 
unexplored. In particular, for the regime $\alpha \in 
[\dims/4, \dims/2]$, the optimal rate of $n^{-1/2}$
for estimating KL divergence (and $D_\phi$ in general) is certainly 
achievable by using Taylor expansion up to second order, 
\emph{assuming} that a separate method exists to achieve 
the optimal rate $n^{-1/2}$ for the integral 
$\int \ptrue^2/\qtrue \psi \;d\mu$.


\section{Conclusions}
\label{Sec:conclude}

We have developed and analyzed $M$-estimation methods for both the
likelihood ratio and $f$-divergence functionals of two unknown
multivariate probability distributions by exploiting a variational
characterization of $f$-divergence functionals. The methods are shown
to be amenable to efficient computational algorithms for optimization
in high-dimensional function spaces. We have also described our method
in the general context of estimating integral functionals of the
likelihood ratio of two unknown densities, and discussed directions
for future work suggested by our results.

\subsection*{Acknowledgements}  This work was partially supported by
NSF grants DMS-0605165 and CCF-0545862 (MJW) and by NSF grant 0509559
(MIJ). The first author would like to acknowledge 
Bharath Sriperumbudur for helpful discussions.


\appendix

\section{Proof of Lemma~\ref{Prop-dist-KL}}
\label{AppLemDistKL}

(i) Note that for $x > 0$, $\frac{1}{2}\log x \leq \sqrt{x} -
1$. Thus, $\int \log \frac{g}{\gtrue} \; d\Pdist \leq 2\int
(g^{1/2}\gtrue^{-1/2} - 1) \; d\Pdist$.  As a result,
\begin{eqnarray*}
d(\gtrue, g) & \geq & \int (g-\gtrue) \; d\Qdist - 
2\int (g^{1/2}\gtrue^{-1/2} - 1) \; d\Pdist \\
& = & \int (g-\gtrue) \; d\Qdist - 
2\int (g^{1/2}\gtrue^{1/2} - \gtrue)\; d\Qdist\\
& = & \int (g^{1/2}-\gtrue^{1/2})^2 d \Qdist.
\end{eqnarray*}
(ii) By our estimation procedure, we have
$\int \gest d\Qdist_n - \int \log \gest d\Pdist_n
\leq \int \gtrue d\Qdist_n - \int \log \gtrue d\Pdist_n$.
It follows that
\begin{eqnarray*}
d(\gtrue, \gest) & = & \int (\gest - \gtrue) d\Qdist
- \int (\log \gest - \log \gtrue) d\Pdist \\
&\leq & \int (\gest - \gtrue) d (\Qdist - \Qdist_n)
- \int (\log \gest - \log \gtrue) d(\Pdist - \Pdist_n) \\
& \leq & \sup_{g\in \Gclass}
\int \log \frac{g}{\gtrue} d(\Pdist_n - \Pdist) -
\int (g - \gtrue) d(\Qdist_n - \Qdist).
\end{eqnarray*} 

\section{Proof of Lemma~\ref{Lem-hellinger}}
\label{AppLemHellinger}

The first inequality is straightforward. We shall focus on the second.
By the definition of our estimator, we have:
\[\int \gest d\Qdist_n - \int \log \gest d \Pdist_n
\leq \int \gtrue d\Qdist_n - \int \log \gtrue d \Pdist_n.\]
Both sides are convex functionals of $g$. Use the following
fact: If $F$ is a convex function and $F(u) \leq F(v)$,
then $F((u+v)/2) \leq F(v)$. We obtain:
\[\int \frac{\gest+\gtrue}{2} d\Qdist_n - \int \log 
\frac{\gest+\gtrue}{2} d \Pdist_n \leq \int \gtrue d\Qdist_n 
- \int \log \gtrue d \Pdist_n.\]
Rearranging,
\begin{eqnarray*}
&& \int \frac{\gest - \gtrue}{2} d(\Qdist_n-\Qdist) -
\int \log \frac{\gest + \gtrue}{2\gtrue} d(\Pdist_n - \Pdist)
 \leq  \int \log \frac{\gest + \gtrue}{2\gtrue} d\Pdist
- \int  \frac{\gest - \gtrue}{2} d\Qdist \\
&& = -d(\gtrue, \frac{\gtrue + \gest}{2}) 
 \leq  - 2h_\Qdist^2(\gtrue, \frac{\gtrue + \gest}{2}),
\end{eqnarray*}
where the last inequality is an application of Lemma~\ref{Prop-dist-KL}.

\comment{ (ii) If $\Gclass$ is a convex class, then $(\gest +
\gtrue)/2 \in \Gclass$.  By the definition of our estimator, there
holds:
\[ \int \gest d\Qdist_n - \int \log \gest d \Pdist_n \leq
\int \frac{\gest+\gtrue}{2} d\Qdist_n - \int \log
\frac{\gest+\gtrue}{2} d \Pdist_n .\] Rearranging,
\[ \int \frac{\gest - \gtrue}{2} d\Qdist_n \leq 
\int \log \frac{2\gest}{\gest + \gtrue} d\Pdist_n
\leq \int \biggr (\frac{2\gest}{\gest + \gtrue} - 1 \biggr ) d\Pdist_n.
\]
It follows that
\begin{eqnarray*}
-\int \frac{\gest - \gtrue}{2} d(\Qdist_n -\Qdist) + \int
\frac{2\gest}{\gest + \gtrue} d(\Pdist_n - \Pdist) & \geq & \int
\biggr (1 - \frac{2\gest}{\gest + \gtrue} \biggr ) d\Pdist + \int
\frac{\gest - \gtrue}{2} d\Qdist \\ = \int \frac{\gtrue - \gest}{\gest
  + \gtrue} d\Pdist + \int \frac{\gest - \gtrue}{2} d\Qdist & = & \int
\frac{\gtrue - \gest}{\gest + \gtrue} \gtrue d\Qdist + \int
\frac{\gest - \gtrue}{2} d\Qdist \\ = \frac{1}{2}\int \frac{(\gtrue
  -\gest)^2}{\gtrue + \gest} d\Qdist & \geq & \frac{1}{2}\int
(\sqrt{\gtrue} -\sqrt{\gest})^2 = h^2 (\gtrue,\gest).
\end{eqnarray*}
}
%
%

\section{Proof of Lemma~\ref{Lem-basic-penalize}}
\label{AppLemBasicPenalize}

Define $d_l(\gtrue,g) = \int (g-\gtrue) d\Qdist - \log\frac{g}{\gtrue}
d\Pdist$.  Note that for $x > 0$, $\frac{1}{2}\log x \leq \sqrt{x} -
1$. Thus,
\begin{eqnarray*}
\int \log \frac{g}{\gtrue} \; d\Pdist & \leq & 2 \int
(g^{1/2}\gtrue^{-1/2} - 1) \; d\Pdist.
\end{eqnarray*}
As a result, for any $g$, $d_l$ is related to $h_\Qdist$ as follows:
\begin{eqnarray*}
d_l(\gtrue, g) & \geq & \int (g-\gtrue) \; d\Qdist - 2\int
(g^{1/2}\gtrue^{-1/2} - 1) \; d\Pdist \\
& = & \int (g-\gtrue) \; d\Qdist - 2\int (g^{1/2}\gtrue^{1/2} -
\gtrue)\; d\Qdist = \int (g^{1/2}-\gtrue^{1/2})^2 d \Qdist \\ & = & 2
h_\Qdist^2(\gtrue,g).
\end{eqnarray*}

By the definition~\eqref{Eqn-penalize} of our estimator, we have:
\[\int \gest d\Qdist_n - \int \log \gest d \Pdist_n 
+ \frac{\lambda_n}{2}I^2(\gest) \leq \int \gtrue d\Qdist_n - \int \log
\gtrue d \Pdist_n + \frac{\lambda_n}{2}I^2(\gtrue).\]
Both sides are convex functionals of $g$. By Jensen's inequality, if
$F$ is a convex function, then $F((u+v)/2) - F(v)\leq (F(u) -
F(v))/2$. We obtain:
\[\int \frac{\gest+\gtrue}{2} d\Qdist_n - \int \log 
\frac{\gest+\gtrue}{2} d \Pdist_n + \frac{\lambda_n}{4}I^2(\gest) \leq
\int \gtrue d\Qdist_n - \int \log \gtrue d \Pdist_n +
\frac{\lambda_n}{4}I^2(\gtrue).\]
Rearranging, $\int \frac{\gest - \gtrue}{2} d(\Qdist_n-\Qdist) - \int
\log \frac{\gest + \gtrue}{2\gtrue} d(\Pdist_n - \Pdist) +
\frac{\lambda_n}{4}I^2(\gest) \leq$
\begin{multline*}
\int \log \frac{\gest + \gtrue}{2\gtrue} d\Pdist - \int \frac{\gest -
  \gtrue}{2} d\Qdist + \frac{\lambda_n}{4}I^2(\gtrue) = -d_l(\gtrue,
  \frac{\gtrue + \gest}{2}) + \frac{\lambda_n}{4}I^2(\gtrue) \\
\leq - 2h_\Qdist^2(\gtrue, \frac{\gtrue + \gest}{2}) +
\frac{\lambda_n}{4}I^2(\gtrue)
\leq - \frac{1}{8} h_\Qdist^2(\gtrue, \gest)
+ \frac{\lambda_n}{4}I^2(\gtrue),
\end{multline*}
where the last inequality is a standard result for the (generalized)
Hellinger distance \\ (cf.~\cite{VdG-book}).

\section{Proof of Theorem~\ref{ThmConvE1}}
\label{AppThmConvE1}

(a) One of the empirical processes on the right-hand side of 
Lemma~\ref{Lem-hellinger} involves function class 
$\Fclass \defn \log \frac{\Gclass + \gtrue}{2\gtrue}$.  
For each $g\in \Gclass$, let $f_g \defn \log
\frac{g+\gtrue}{2\gtrue}$.  We endow $\Fclass$ with a ``norm,''
namely, \emph{Bernstein distance}.  This is defined as
follows: for a constant $K > 0$,
\[\rho_K(f)^2 \defn 2K^2\int (e^{|f|/K} - 1 - |f|/K) d\Pdist.\]

The Bernstein distance is related to the Hellinger distance in several
crucial ways (see, e.g.,~\cite{VdG-book}, page 97):
\begin{itemize}
\item $\rho_1(f_g) \leq 4h_\Qdist(\gtrue,\frac{g+\gtrue}{2})$.
\item The bracketing entropy based on Bernstein distance is also
related to the bracketing entropy based Hellinger distance
(i.e., which is the $L_2$ norm for the square root function):
\begin{equation}
\label{Eqn-bernstein-entropy}
\entrobrack{\sqrt{2}\delta}{\Fclass}{\rho_1}
\leq \entrobrack{\delta}{\bar{\Gclass}}{L_2(\Qdist)},
\end{equation}
where $\bar{\Gclass} \defn \{((g+\gtrue)/2)^{1/2}, g\in \Gclass\}$
and $\bar{g} \defn (g+\gtrue)/2$.
\end{itemize}

By Lemma~\ref{Lem-hellinger}, for any $\delta > 0$, with respect to 
$\Pdist$ measure:
\begin{eqnarray*}
\lefteqn{P(h_\Qdist(\gtrue,\gest) > \delta) \leq 
P(h_\Qdist(\gtrue,(\gest + \gtrue)/2) > \delta/4) } \\
& \leq & P \biggr ( \sup_{g\in \Gclass,\; h_\Qdist(\gtrue,\bar{g}) 
> \delta/4} -\int (\bar{g} - \gtrue) d(\Qdist_n - \Qdist)
+ \int f_g \;d(\Pdist_n - \Pdist) - 2h_\Qdist^2(\gtrue,\bar{g}) 
\geq 0 \biggr ) \\
& \leq & P \biggr (\sup_{g\in \Gclass,\; h_\Qdist(\gtrue,\bar{g}) 
> \delta/4} -\int (\bar{g} - \gtrue) d(\Qdist_n - \Qdist) - 
h_\Qdist^2(\gtrue,\bar{g}) \geq 0 \biggr) \\
&& + \; P \biggr (\sup_{g\in \Gclass,\; h_\Qdist(\gtrue,\bar{g}) 
> \delta/4} \int f_g \;d(\Pdist_n - \Pdist) - 
h_\Qdist^2(\gtrue,\bar{g}) \geq 0 \biggr )
\defn A + B.
\end{eqnarray*}
We need to upper bound the two quantities $A$ and $B$ on the
right-hand side of this equation.  These can be handled in a similar manner.
Since $\entrobrack{\delta}{\bar{\Gclass}}{L_2(\Qdist)} < \infty$ 
the diameter of $\bar{\Gclass}$ is finite.
Let $S$ be the minimum $s$ such that $2^{s+1}\delta/4$ exceeds that diameter.
We apply the so-called peeling device: Decompose
$\bar{\Gclass}$ into layers of Hellinger balls around $\gtrue$ and then
apply the union bound on the probability of excess. For each layer,
one can now apply the modulus of continuity of suprema of an
empirical process.
\begin{eqnarray*}
B & \leq & \sum_{s=0}^{S} 
P \biggr (\sup_{g\in \Gclass,\; h_\Qdist(\gtrue,\bar{g}) \leq 2^{s+1}\delta/4}
\int f_g \;d(\Pdist_n - \Pdist) \geq 2^{2s}(\delta/4)^2 \biggr ).
\end{eqnarray*}
Note that if $h_\Qdist(\gtrue,\bar{g}) \leq 2^{s+1}\delta/4$
then $\rho_1(f_g) \leq 2^{s+1}\delta$. Note that for any $s=1,\ldots,S$,
the bracketing entropy integral can be bounded as:
\begin{eqnarray*}
\lefteqn{\int_{0}^{2^{s+1}\delta}\entrobrack{\epsilon}{\Fclass \cap 
\{h_\Qdist(\gtrue,\bar{g}) \leq 2^{s+1}\delta/4\}}{\rho_1}^{1/2} \;d\epsilon} \\
& \leq & 
\int_{0}^{2^{s+1}\delta}\entrobrack{\epsilon/\sqrt{2}}{\bar{\Gclass} \cap
\{h_\Qdist(\gtrue,\bar{g}) \leq 2^{s+1}\delta/4\}
}{L_2(\Qdist)}^{1/2} \;d\epsilon \\
& \leq & \int_{0}^{2^{s+1}\delta} 
C_9(\epsilon/\sqrt{2})^{-\gamma_{\bar{\Gclass}}/2} 
\; d\epsilon \\
& \leq & C_8 (2^{s+1}\delta)^{1-\gamma_{\bar{\Gclass}}/2},
\end{eqnarray*}
where $C_8,C_9$ are constants independent of $s$. Now apply
Theorem~\ref{Thm-VdG-rate} (see Appendix~\ref{AppEmpProcess}), where
$K=1$, $R = 2^{s+1}\delta$, $a = C_1 \sqrt{n}R^2/K = C_1
\sqrt{n}2^{2(s+1)}\delta^2$. We need
\[a \geq C_0 C_8(2^{s+1}\delta)^{1-\gamma_{\bar{\Gclass}}/2} > C_0 R.\]
This is satisfied if $\delta = n^{-1/(\gamma_{\bar{\Gclass}} + 2)}$
and $C_1 = C_0 C_8$, where $C_8$ is sufficiently large (independently
of $s$).  Finally, $C_0^2 \geq C^2(C_1 + 1) = C^2(C_0 C_8 + 1)$ if
$C_0 \defn 2 C^2 C_8 \vee 2C$, where $C$ is some universal constant in
Theorem~\ref{Thm-VdG-rate}. Applying this theorem, we obtain:
\begin{eqnarray*}
B & \leq & \sum_{s=0}^{S} C \exp \biggr [-\frac{C_1^2 n 2^{2(s+1)}\delta^2}
{C^2(C_1 + 1)} \biggr] \leq c\exp\biggr [-\frac{n\delta^2}{c^2} \biggr ]
\end{eqnarray*}
for some universal constant $c$.  A similar bound can be obtained for 
$A$, with respect to $\Qdist$ measure and with $\delta =
n^{-1/(\gamma_{\bar{\Gclass}}+2)}$.
Since $\ptrue/\qtrue$ is bounded from above, this also implies a
probability statement with respect to $\Pdist$.  Thus,
$h_\Qdist(\gtrue,\gest)$ is bounded in $\Pdist$-probability by
$n^{-1/(\gamma_{\bar{\Gclass}}+2)}$.  \\ (b) The proof is similar to
Theorem~\ref{ThmConvE2}(b) and is omitted.

\section{Comparison of the rate in Lemma~\ref{Prop-dist-KL} to the minimax rate}
\label{AppMiniMax}

Recall that the minimax rate is defined as
\begin{equation*}
r_n \; \defn \; \inf_{\gest \in \Gclass} \sup_{\Pdist,\Qdist}
\E_\Pdist\, [h_\Qdist(\gtrue, \gest)],
\end{equation*}
where the supremum is taken over all pairs $(\Pdist, \Qdist)$ such
that $\gtrue \in \Gclass$.  Note that $r_n \geq \inf_{\gest \in \Gclass}
\sup_{\Pdist} \E h_\mu(\gtrue,\gest)$, where we have fixed $\Qdist =
\mu$, the Lebesgue measure on $\Xspace$.  We can reduce this 
bound to the minimax lower bound for a nonparametric density
estimation problem~\citep{Yu}.  This reduction is not entirely
straightforward, however, because the space $\Gclass$ ranges 
over smooth functions that need not be valid probability 
densities.  Therefore, an easy-to-use minimax lower bound 
such as that of~\citep{Yang-Barron} is not immediately
applicable. Nonetheless, we can still apply the hypercube argument and
the Assouad lemma to obtain the right minimax rate. See \S 24.3 of van
der Vaart~\cite{vanderVaart98} for a proof for the case of one
dimension. This proof goes through for general $\dims \geq 1$.


\section{Some calculations for Theorem~\ref{ThmExt}}
\label{AppCalc}

Note that the conjugate dual of $\phi(u) = 1/u$ takes the form
\begin{eqnarray*}
\phi^*(v) & = & \begin{cases} -2\sqrt{-v} & \mbox{if $v< 0$, and} \\
 +\infty & \mbox{otherwise.}
\end{cases}
\end{eqnarray*}
Consequently, we can restrict $\Fclass$ to the subset for which $f <
0$ for any $f\in \Fclass$.  Let $g \defn \sqrt{-f}$ and $\Gclass =
\sqrt{-\Fclass}$.  $\Gclass$ is a function class of positive
functions.  We have $\gtrue \defn \sqrt{-\ftrue} =
\sqrt{-\phi'(\qtrue/\ptrue)} = \ptrue/\qtrue$. Define $\gest \defn
\sqrt{-\fest}$.  We also replace notation $d_\phi(\ftrue,f)$ by
$d_\phi(\gtrue,g)$.  For our choice of $\phi$, we have:
\begin{eqnarray*}
d_\chi(\gtrue,g) & = & d_\chi(\ftrue, f) = \int (-2\sqrt{-f} +
2\sqrt{-\ftrue}) d\Pdist -(f-\ftrue)d\Qdist \\ & = & \int (\gtrue -
g)(2\ptrue/\qtrue - \gtrue - g) d\Qdist \\ & = & \int(g - \gtrue)^2
d\Qdist,
\end{eqnarray*}
as claimed.
Moreover, we have
\begin{eqnarray*}
v_n^\chi(\Gclass) & = & v_n^\chi(\Fclass) = \sup_{g \in \Gclass}
\biggr | \int 2(g^2 - \gtrue^2) d(\Qdist_n - \Qdist) - \int (g -
\gtrue) d(\Pdist_n - \Pdist) \biggr |.
\end{eqnarray*}

\section{Results from empirical process theory}
\label{AppEmpProcess}

For completeness, we state here two standard results from
empirical process theory that are needed in the paper.  These
results are versions of Theorems 3.7 and 5.11 from van de
Geer~\cite{VdG-book}, respectively:
\begin{theorem}
\label{Thm-VdG-consistency}
Let $G(x) = \sup_{g \in \Gclass} |g(x)|$ be the envelope function for
a function $\Gclass$.  Assume that $\int G d\Pdist < \infty$, and
suppose moreover that for any $\delta > 0$,
$\frac{1}{n}\entropy{\delta}{\Gclass}{L_1(\Pdist_n)} \convp{\Pdist}
0$.  Then $\sup_{g\in \Gclass} \int g d(\Pdist_n - \Pdist) \convas 0$.
\end{theorem}

\begin{theorem}
\label{Thm-VdG-rate}
Suppose that the function class $\Gclass$ satisfies $\sup_{g\in
\Gclass}\rho_K(g) \leq R$ for some constants $K$ and $R$.  Given $a > 0$,
suppose that for some constants $C$ and $C_1$, there holds
\begin{eqnarray*}
a & \leq & C_1 \sqrt{n}R^2/K \\ a & \geq & \sqrt{ C^2(C_1 + 1)} \;
\biggr (\int_{0}^{R}\entrobrack{u}{\Gclass}{\rho_K}^{1/2}du \vee R
\biggr ).
\end{eqnarray*}
Then the empirical process is bounded as
\begin{eqnarray}
\Prob \biggr[\sup_{g\in \Gclass} \big|\sqrt{n}\int g d(\Pdist_n -
 \Pdist) \big| \geq a \biggr ] & \leq & C \exp \biggr
 [-\frac{a^2}{C^2(C_1 + 1)R^2} \biggr ].
\end{eqnarray}
\end{theorem}


\bibliography{mutin}


\end{document}